\documentclass[twoside]{article}
\usepackage{latexsym, amssymb, enumerate, amsmath, amsthm}
\usepackage{latexsym,amssymb}

\usepackage[top=2cm, bottom=2cm, left=3cm, right=3cm]{geometry}

\newcommand{\fer}[1]{(\ref{#1})}

\newcommand{\norm}[1]{\|#1 \|}
\newcommand{\R}{\mathbb{R}}

\newcommand{\Sph}{\mathbb{S}}

\newcommand {\e}  {\varepsilon}
\newcommand {\eps}  {\varepsilon}

\newcommand {\cad} { {\cal D} }

\newcommand {\f}   {\frac}

\newcommand{\beq}{\begin{equation}}
\newcommand{\beqa}{\begin{eqnarray}}
\newcommand{\bea} {\begin{array}{ll}}
\newcommand{\beqan}{\begin{eqnarray*}}
\newcommand{\eeq}{\end{equation}}
\newcommand{\eeqa}{\end{eqnarray}}
\newcommand{\eeqan}{\end{eqnarray*}}
\newcommand{\eea} {\end{array}}

\newcommand{\bone}{\mathbf{1}}

 \numberwithin{equation}{section}

\newtheorem{thm}{Theorem}[section]
\newtheorem{proposition}[thm]{Proposition}
\newtheorem{lem}[thm]{Lemma}

\newtheorem{rem}[thm]{Remark}

\begin{document}

\title{The parabolic-parabolic Keller-Segel model in $\R^2$
\footnote{The authors warmly thank Beno\^it Perthame for challenging discussions. The authors  acknowledge the support from European Union Marie Curie network M3CS-TuTh, MRTN-CT-2004-503661, ANR project MACBAC.
}}
\author{Vincent Calvez
\footnote{D\'epartement de Math\'ematiques et Applications, 
\'Ecole Normale Sup\'erieure, CNRS UMR8553 , 
            45 rue d'Ulm, F~75230 Paris cedex 05 
(email: vincent.calvez@ens.fr).}
and Lucilla Corrias \footnote{D\'epartement de Math\'ematiques, 
Universit\'e d'Evry Val d'Essonne, 
Rue du P\`ere Jarlan, F~91025 Evry Cedex 
(email: lucilla.corrias@univ-evry.fr).}}

\date{}

\pagestyle{myheadings} 
\markboth{The parabolic-parabolic Keller-Segel model in $\R^2$}
{V. Calvez and L. Corrias}
\maketitle

\begin{abstract}
This paper is devoted mainly to the global existence problem for the two-dimensional parabolic-parabolic Keller-Segel in the full space. We derive a critical mass threshold below which global existence is ensured. Using carefully energy methods and {\em ad hoc} functional inequalities we improve and extend previous results in this direction. The given threshold is supposed to be the optimal criterion, but this question is still open. This global existence result is accompanied by a detailed discussion on the duality between the Onofri and the logarithmic Hardy-Littlewood-Sobolev inequalities that underlie the following approach. \end{abstract}

\paragraph{keywords.} Chemotaxis, parabolic system, global weak solutions, energy method, Onofri inequality, Hardy-Littlewood-Sobolev inequality.
\smallskip


\section{Introduction}
\label{Sec:intro}

Within living organisms, cells may communicate and therefore interact through chemical signals. This signaling pathway is of crucial importance for cell particles to move in the right direction or to organize themselves spatially. Biological challenges involving this phenomenon are numerous. It is known to be so actually in immunology and inflammatory processes, in bacterial growth colony and at some key stages of embryonic development for instance. Among mathematical models describing spatial organization of biological population through chemical signals (\cite{P}, \cite{Markowich}, \cite{ErbanOthmer}) we highlight the following Patlak-Keller-Segel (PKS) model for chemotaxis (\cite {KS70},\cite{Patlak53}) 
\begin{equation}
\left\{\begin{array}{rcl}
\dfrac{\partial n}{\partial t}  & =&\kappa \Delta n - \chi\nabla\cdot( n\nabla c)\ ,\quad t>0,\ x\in \Omega,\vspace{.2cm}\\
\e\dfrac{\partial c}{\partial t} &=& \eta\Delta c + \beta n - \alpha c\ ,\quad t>0,\ x\in \Omega,
\end{array}
\right.  
\tag{PKS}
\end{equation}
where $n$ denotes the density of a cell population and $c$ is the concentration of a chemical signal attracting the cells. The parameters $\kappa$, $\chi$ and $\beta$ are given positive constants, while $\e$, $\eta$ and $\alpha$ are given non-negative constants determining the type of evolution undergone by $c$. The set $\Omega$ is either a bounded domain in $\R^d$ or the whole space $\R^d$. In any case, boundary conditions or decay conditions at infinity have to be given for the densities $n$ and $c$ together with the initial conditions $n(\cdot,0)=n_0$ and $c(\cdot,0)=c_0$ if $\e>0$. 

The modeling interest of such a system is to exhibit a phenomenon of ``critical mass'' at least in dimension $d=2$. In a system like (PKS), the coupling between the cell equation and the chemical equation is a positive feedback: the more cells are aggregated, the more they produce a signal attracting other cells. This process is counter-balanced by pure diffusion of the cells, but if the amount of cells is sufficient this non-local chemical interaction dominates and cells attract themselves. Thus model (PKS) provides a simple phenomenological description of an instability broadly encountered in biology. The starvation stage of the slime mold amoebae {\em Dictyostelium discoideum} for instance is governed by this process, driving the population of unicellular organisms into a multicellular one (\cite{Weijer}, \cite{Hofer}). More recently, a (PKS) type model has been suggested to solve a remarkable pattern formation issue in the human brain \cite{KC07}. 

It is not overemphasize to say that system (PKS) has been the subject of a huge quantity of mathematical analysis over the last thirty years. The results of all these investigations can be simply summarized saying that the global existence or the blow-up of solutions of (PKS) is a space dimension dependent phenomenon. In particular, in dimension $d=2$, the above biological instability has been precisely described mathematically at least when  $\e=\alpha=0$. Indeed, for the parabolic-elliptic system it has been shown that there exists a threshold for the initial mass $M=\int_{\Omega}n_0(x)\ dx$. For values of $M$ under this threshold the solution exists globally in time, while above this threshold the solution blow-up in finite time (\cite{B98},\cite{BDP},\cite{Gajewski98},\cite{NSY}). The modality of the blow-up has also been analyzed \cite{HerreroVelazquez}, as well as the critical case, i.e. when the initial mass $M$ equals the threshold value (\cite{BCM}, \cite{Biler06}).  

Let us mention that in dimension $d\ge3$ a similar critical phenomenon has been investigated. In this case, the $L^{\f d2}$-norm of the initial density $n_0$ plays the same role as the initial mass $M$ in dimension $d=2$. Indeed, in \cite{CP_CRAS06}, \cite{CP07}, \cite{CPZ04} the authors proved the global existence of weak solution of system (PKS) of parabolic-parabolic, parabolic-elliptic and parabolic degenerate type under a smallness condition on $\| n_0\|_{L^{d/2}}$. However, up to our knowledge, blow-up for large $\| n_0\|_{L^{d/2}}$ is still open and we conjecture that no critical threshold exists as in dimension 2.

Despite of all these results, mathematical open problems around (PKS) still subsist, especially  for the full parabolic-parabolic system (PKS) ($\e>0$). Let us observe that, whenever $\e=0$ (quasi-stationary hypothesis for the chemical $c$), the system reduces to a single parabolic equation with a quadratic nonlocal nonlinearity, $c$ being expressed as a convolution between $n$ and a kernel.  On the other hand, when $\e>0$ the full parabolic-parabolic system is more difficult to handle with, since the equations are strongly coupled. However, we will see that in this case the densities $n$ and $c$ play dual roles in some sense (as it is highlighted by the dual inequalities used throughout this paper), providing some interesting features and structure to the system (PKS).

The goal of this paper is to tackle the global existence problem for the full parabolic-parabolic system (PKS) with $\e>0$, $\eta>0$, $\alpha\ge0$ in dimension $d=2$ and in the whole space $\R^2$. Indeed, for this problem the optimal threshold of $M$ for the global existence of solutions has not been found yet. A result exists in this direction in \cite{N01}, but it does not give the exact critical mass. Here we obtain the optimal critical mass value using the energy method (\cite{Gajewski98},\cite{B98},\cite{BDP}) and {\em ad hoc} functional inequalities on $\R^2$. One more time, the free energy functional 
\begin{equation*}
\begin{split}
\mathcal E(t) = \int_{\R^2} n(x,t)\log n(x,t)\ dx - \int_{\R^2} n(x,t)c(x,t)\ dx \\ + \f12\int_{\R^2} |\nabla c(x,t)|^2\ dx + \f\alpha2 \int_{\R^2} c^2(x,t)\ dx \label{eq:free energy}
\end{split}
 \end{equation*}
comes out to be the key ingredient leading to the  global   existence of solutions under the optimal smallness condition for the mass. Indeed, $\mathcal E(t)$ together with its evolution equation provide a gallery of {\em a priori} estimates on the solutions $(n,c)$ and we shall make use of each of them. For instance, they allow to prove that the cellular flux in \eqref{KS} $n\big(\nabla(\log n-c)\big)\in L^1(\R^+_{loc}\times\R^2)$. Therefore, the equation on $n$ holds in the distribution sense. Surprisingly, no specific restrictions on $c_0$ are required even for the fully parabolic-parabolic case under interest, except of course suitable regularity of the initial data.

It is convenient to adimensionalize system (PKS) through the following change of variables
$$
t\to\tau= \kappa t\ ,\quad n(x,t)\to\tilde n(x,\tau)=\f{\beta\chi}{\eta\kappa}\ n\left(x,\f\tau\kappa\right)\ ,\quad
c(x,t)\to\tilde c(x,\tau)=\f{\chi}{\kappa}\ c\left(x,\f\tau\kappa\right)\ .
$$
Therefore, the system under consideration will be
\begin{equation}
\left\{\begin{array}{rcl}
\partial_t \tilde n& =& \Delta \tilde n - \nabla\cdot(\tilde n\nabla\tilde c)\ ,\quad t>0,\ x\in \R^2,\vspace{.2cm}\\
\tilde\e\partial_t\tilde c&= &\Delta\tilde c + \tilde n - \tilde\alpha\tilde c\ , \quad\quad t>0,\ x\in \R^2,\vspace{.2cm}\\
\tilde n(\cdot,0)&=&\tilde n_0(\cdot)=\f{\beta\chi}{\eta\kappa}\ n_0(\cdot)\ ,\quad \tilde c(\cdot,0)=\tilde c_0(\cdot)=\f{\chi}{\kappa}\ c_0(\cdot),\quad x\in \R^2\ .
\end{array}
\right. 
\label{KS}
\end{equation}
with $\tilde\e=\f{\e\kappa}\eta>0$ and $\tilde\alpha=\f\alpha\eta\ge0$. The tilde sign will be removed in the sequel for clarity. Moreover, fast decay conditions at infinity for $n$ and $c$ have to be associated with \fer{KS}. Concerning the cell density, this decay will be expressed in terms of moments of $n$.

After this change of variables, the only parameters of the system to deal with  are the total mass of cells $M=\int_{\R^2}n_0(x)\ dx$, which is conserved along time, the inverse  diffusion rate of the chemical $\e$ and the chemical degradation rate $\alpha\geq0$. The latter seems to play no essential role, unless it induces slightly technical difficulties in the estimates. 

Our main results are the followings.
\begin{thm}[Global existence]
\label{the:existence}
Assume $\e>0$ and $\alpha\geq 0$. Let $(n_0,c_0)$ be non-negative initial conditions for the parabolic-parabolic system \fer{KS} such that 
\begin{description}
\item[(H1)] $n_0\in L^1(\R^2)\cap L^1(\R^2,\log(1+|x|^2)dx)$ and $n_0\log n_0\in L^1(\R^2)$\ ; 
\item[(H2)] $c_0\in H^1(\R^2)$ if $\alpha>0$ while $c_0\in L^1(\R^2)$ and $|\nabla c_0|\in L^2(\R^2)$  if $\alpha=0$\ ;
\item[(H3)] $n_0c_0\in L^1(\R^2)$\ .
\end{description}
\noindent Assume in addition that the mass is subcritical, i.e. $M<8\pi$. Then, there exists a global weak non-negative solution $(n,c)$ of \fer{KS} such that 
\[
\begin{split}
n\in L^\infty((0,\infty);L^1(\R^2))\cap L^\infty_{loc}((0,\infty);L^1(\R^2,\log(1+|x|^2)dx))\\ \hbox{and}\quad
n\log n\in L^\infty_{loc}((0,\infty);L^1(\R^2))\ ;
\end{split}
\]
$$
c\in L^\infty_{loc}((0,\infty);H^1(\R^2))\quad\hbox{if }\alpha>0\ ;
$$
$$
c\in L^\infty_{loc}((0,\infty);L^1(\R^2))\quad\hbox{and}\quad 
|\nabla c|\in L^\infty_{loc}((0,\infty);L^2(\R^2))
\quad\hbox{if }\alpha=0\ ;
$$
$$
\partial_tc\in L^2_{loc}((0,\infty);L^2(\R^2))\quad\hbox{and}\quad
nc\in L^\infty_{loc}((0,\infty);L^1(\R^2)) \ ;
$$
$$
\int_0^t\int_{\R^2} n(x,s)|\nabla(\log n(x,s) - c(x,s)) |^2 \ dxds<\infty\quad\forall t>0\ ;
$$
$$
\mathcal E(t)+ \int_{0}^{t}\! \int_{\R^2} n(x,s)|\nabla(\log n(x,s) - c(x,s)) |^2 \  dxds
+\e\int_{0}^{t}\!\int_{\R^2} |\partial_t c(x,s)|^2\ dxds  \le\mathcal E(0)\ .
$$
Moreover, $n\in L^\infty_{loc}((0,\infty);L^p(\R^2))$ for any $1<p<\infty$ (regularizing effect).
\end{thm}\
%

%
Let us observe that the global existence result given in Theorem \ref{the:existence} holds true under the optimal  condition $n_0\in L^1(\R^2,\log(1+|x|^2)dx)$. This is weaker than the hypothesis $n_0\in L^1(\R^2,|x|^2dx)$ used in  \cite{BDP} for the global existence result in the case $\e=\alpha=0$. As a matter of fact, we can adapt Theorem~\ref{the:existence} to the case $\e=\alpha=0$ and slightly improve the result in \cite{BDP}, (see Appendix \ref{ap:eps=alpha=0}). On the other hand, the hypothesis $n_0\in L^1(\R^2,|x|^2dx)$ remains still necessary to obtain a blow-up result for the parabolic-elliptic system \fer{KS}. Nevertheless, we can extend the blow-up result in \cite{BDP} to the case  $\alpha\ge0$ as follows. The same blow-up result for the full parabolic-parabolic system is still, at our knowledge, an open problem.


%
\begin{thm}[Blow-up]
Assume $\eps=0$ and $\alpha\ge0$. Let $n_0$ be a non-negative $L^1(\R^2)$ function with super-critical mass, i.e. $M>8\pi$ and finite second momentum $I_0$. Let $n$ be a non-negative smooth solution of \fer{KS} with $c=B_\alpha*n$ and $B_\alpha$ the Bessel kernel, and let $[0,T^*)$ be the maximal interval of existence. Then, there exists a positive universal constant $\mathcal C$ such that if 
\begin{equation}
\alpha \int_{\R^2} |x|^2n_0(x)\ dx \leq \f1{4\mathcal C^2 M}(M-8\pi)^2 \ ,  
\label{I(0)}
\end{equation}
the solution blow-up and $T^*\le\f{2\pi I_0}{M(M-8\pi-2\mathcal C\sqrt{\alpha MI_0})}$\ .
\label{the:blowup}
\end{thm}

The paper is organized as follows. In Section \ref{Sec:free energy} we give a set of technical tools to be used in the sequel and we describe briefly the dual strategies that give the required {\em a priori} estimates. In Section~\ref{sec:Onofri} we derive the {\em a priori} estimates in both cases $\alpha>0$ and $\alpha=0$ (amoung wich the key equi-integrability estimate) from the so-called Onofri's inequality on the whole space $\R^2$ . In Section~\ref{sec:log HLS} we re-derive those estimates using the dual strategy based on the logarithmic Hardy-Littlewood-Sobolev inequality. Section \ref{Sec:regularizing} is devoted to the proof of the regularizing effect acting on the solutions. Section~\ref{Sec:Globalexistence} is a short description of the regularization procedure which leads to the rigorous proof of global existence when combined with the {\em a priori} estimates derived in Sections \ref{sec:Onofri} and \ref{sec:log HLS}. Blow-up of the solutions in the special case $\e=0$ is shown for a super-critical mass in Section \ref{Sec:blowup}, under a smallness condition on $\alpha\int_{\R^2}|x|^2n_0(x)\ dx$. Finally, several complementary results are given for the sake of completeness in Section \ref{sec:appendix}.

In the sequel, we will denote by $C$ every positive constants that may vary from line to line in the computations. Only the dependence on crucial parameters will be written explicitly.  Moreover, $(u)_+$ and $(u)_-$ will denote the positive and negative part of $u$ as usual, while $L^1_+(\R^2):=\{f\in L^1(\R^2)\ ,\ f\ge0\}$.

\section{The free energy and the moments control}
\label{Sec:free energy}
It is well known that system \fer{KS} is equipped with the following {\sl free energy functional}
\begin{equation} \begin{split}
\mathcal E(t) = \int_{\R^2} n(x,t)\log n(x,t)\ dx - \int_{\R^2} n(x,t)c(x,t)\ dx \\ + \f12\int_{\R^2} |\nabla c(x,t)|^2\ dx + \f\alpha2 \int_{\R^2} c^2(x,t)\ dx\ .
\label{eq:fe}
\end{split}
\end{equation}
In the kinetic equation literature, the first term $\int_{\R^2} n(x,t)\log n(x,t)\ dx$ is usually referred to as the physical entropy. However, here it will be more convenient and natural to define the {\sl entropy} in the line of \cite{CJMTU} as
\begin{equation}
E(n;c)(t)=\int_{\R^2} (n(x,t)\log n(x,t)-n(x,t)c(x,t))\ dx\ ,
\label{entropy}
\end{equation}
including also the potential energy term $\int_{\R^2} n(x,t)c(x,t)dx$. On the other hand, the potential energy term has also to be included in the {\sl chemical energy} associated to the elliptic equation $-\Delta c+\alpha c=n$, i.e.
\begin{equation}
F_\alpha(c;n)(t)=\f12\int_{\R^2} |\nabla c(x,t)|^2\ dx + \f\alpha2 \int_{\R^2} c^2(x,t)\ dx-\int_{\R^2} n(x,t)c(x,t)\ dx\ ,\quad\alpha\ge0\ .
\label{chemicalenergy}
\end{equation}
Thereby, the free energy $\mathcal E(t)$ is a superposition of \fer{entropy}  and \fer{chemicalenergy} thus reflecting the strongly coupled property of system \fer{KS}. The quantity $\mathcal E(t)$ will play a fondamental role in the research of {\em a priori} estimates starting from the following proposition.

\begin{proposition}\label{prop:fe}
Let $(n,c)$ be any non-negative and sufficiently smooth solution of \fer{KS} with finite free energy \fer{eq:fe}. Then $\mathcal E(t)$ decreases along the trajectories of the dynamical system associated to \fer{KS}, since 
\begin{equation}
\dfrac{d}{dt}\mathcal E(t) = - \int_{\R^2} n(x,t)|\nabla(\log n(x,t) - c(x,t)) |^2 \ dx - \e\int_{\R^2} |\partial_t c(x,t)|^2\ dx\leq 0\ .
\label{eq:energy dissipation}
\end{equation}
\end{proposition}
\begin{proof}
The equation on $n$ can be written as $\partial_t n=\nabla\cdot\big(n\nabla(\log n-c)\big)$. Then, using the mass conservation, we obtain
\begin{eqnarray} 
\int_{\R^2} \partial_t n(x,t)\ (\log n(x,t) -c(x,t))\ dx& =& 
\dfrac d {dt} \int_{\R^2} n(x,t)\log n(x,t)\ dx - \int_{\R^2}  \partial_t n(x,t) \ c(x,t)\ dx \nonumber\\
&=&  - \int_{\R^2} n(x,t)|\nabla (\log n(x,t) -c(x,t))|^2\ dx\ .
\label{eq:1}
\end{eqnarray}
On the other hand, testing the equation on $c$ against $\partial_t c$, we have
\begin{equation}
\begin{split}
\e\int_{\R^2} |\partial_t c(x,t)|^2\ dx = -\frac d {dt}\int_{\R^2} \f{|\nabla c(x,t)|^2}2\ dx + \int_{\R^2} n(x,t)\ \partial_t c(x,t)\ dx\\ -\alpha \frac d{dt} \int_{\R^2} \f{c^2(x,t)}2\ dx\ . 
\label{eq:2}  
\end{split}
\end{equation}
We conclude by summing \fer{eq:1} and \fer{eq:2}.
\end{proof}

Equation \fer{eq:energy dissipation} measures the dissipation of the free energy due to the {\sl entropy production} term 
\begin{equation}
I(t)=\int_{\R^2} n(x,t)|\nabla(\log n(x,t) - c(x,t)) |^2 \ dx \ ,
\label{entropyproduction}  
\end{equation}
and to the chemical production term $\e\int_{\R^2} |\partial_t c(x,t)|^2\ dx$.  Let us observe however that any weak solution of \fer{KS} is not expected to satisfy \fer{eq:energy dissipation} but the inequality
\begin{equation}
\mathcal E(t)+ \int_{0}^{t} \int_{\R^2} n(x,s) |\nabla(\log n(x,s) - c(x,s)) |^2 \  dx  ds
+\e\int_{0}^{t}\int_{\R^2} |\partial_t c(x,s)|^2\ dx  ds  \le\mathcal E(0)\ ,
\end{equation}
as it is under the quasi-stationary hypothesis $\e=\alpha=0$ (see \cite{BDP}).

The time-monotonicity of $\mathcal E(t)$ given by \eqref{eq:energy dissipation} provides us with an upper control of the entropy \fer{entropy}. But a control from below of the entropy is also needed in order to obtain {\em a priori} estimates on the solution and then the global existence result. In the case of the parabolic-parabolic system \fer{KS} on a bounded domain $\Omega\subset\R^2$, the strategy usually followed makes use essentially of two primary tools: a minimization principle with respect to $n$ of the entropy $E(n;c)$ and the Moser-Trudinger inequality. Moreover, their combination gives the {\em a priori} estimates under the exact critical mass value for $M$ (see \cite{B98}, \cite{calvezc}, \cite{Gajewski98}, \cite{NSY}).

Concerning system \fer{KS} in the whole space $\R^2$ under interest, the first result we give here is that the same method above can be followed. However, in order to do so, one has firstly to reinforce the space decay of $c$ as $|x|\to+\infty$ in order to minimize the entropy $E(n;c)$ with respect to $n$ and secondly to employ  an {\em ad hoc} Moser-Trudinger type inequality, i.e. the Onofri's inequality \cite{Ono}. 

The second result is that an alternative strategy, dual in some sense to the previous one, can be also adopted. The technical tools to be employed are: the minimization with respect to $c$ of the chemical energy $F_\alpha(c;n)$ \fer{chemicalenergy} and the logarithmic Hardy-Littlewood-Sobolev inequality (HLS in the sequel, see Lemma \ref{lem:logHLS} and  \cite{CarlenLoss}) if $\alpha=0$, or a modified version of this inequality if $\alpha>0$ (see Lemma \ref{lem:Besselinequality}). This method is new and it is somewhat the extension to the parabolic-parabolic system \fer{KS} of what was done in \cite{BDP} for the parabolic-elliptic system  \fer{KS} with $\e=\alpha=0$. Indeed, in this case the free energy reads as
\begin{equation}
\mathcal E(t) = \int_{\R^2} n(x,t)\log n(x,t)\ dx-\f12\int_{\R^2} n(x,t)c(x,t)\ dx\ ,
\label{freeellipticenergy}
\end{equation}
with the concentration of the chemical given by $c(x,t)=-\f1{2\pi} \int_{\R^2}\log|x-y|n(y,t)\ dy$, and hence it is well adapted to apply  the HLS inequality. 

No matter of the method followed to obtain the necessary {\em a priori} estimates, we are in any case lead to consider the following modified free energy functional 
\begin{eqnarray*}
\mathcal E_H(t) &=&\mathcal E(t)-\int_{\R^2}n(x,t)\log H(x)\ dx \nonumber \\ 
&=&E(n;c + \log H)+ \f12\int_{\R^2} |\nabla c(x,t)|^2\ dx + \f\alpha2 \int_{\R^2} c^2(x,t)\ dx\ ,
\label{eq:feH}
\end{eqnarray*}
where 
\begin{equation}
H(x)=\f1\pi\f1{(1+|x|^2)^2}
\label{H}  
\end{equation}
has been chosen so that ${\cal J_S}:=4\pi H$ is the Jacobian of the usual stereographic projection on the sphere ${\cal S}\ :\ \R^2\cup\{\infty\}\to\Sph^2$ and $\int_{\R^2}H(x)dx=1$, (see \cite{LiebLoss}).

The introduction of the function $H$ will appear to the reader more natural in Sections \ref{sec:Onofri} and  \ref{sec:log HLS} where the two methods will be developed respectively. Here, let us observe that, by opposition to $\mathcal E(t)$, the functional $\mathcal E_H(t)$ is not time decreasing.  However, we can control its time-growth by the following computation. We have 
\begin{eqnarray*}
 \f{d}{dt}  \int_{\R^2}\ n(x,t)\log H(x)\ dx &=&-\int_{\R^2} n(x,t)\nabla \log H(x)\cdot \nabla (\log n(x,t) - c(x,t))\ dx\\
 &=&2\int_{\R^2} n(x,t)\nabla \log(1+|x|^2)\cdot \nabla (\log n(x,t) - c(x,t))\ dx 
 \end{eqnarray*}
and using equation \fer{eq:energy dissipation}, we easily obtain
\begin{eqnarray*}
\dfrac{d}{dt}\mathcal E_H(t) &=& - \int_{\R^2} n(x,t)|\nabla(\log n(x,t) - c(x,t)) |^2 \ dx - \e\int_{\R^2} |\partial_t c(x,t)|^2\ dx\\  
& &\quad\quad-2\int_{\R^2} n(x,t)\nabla \log(1+|x|^2)\cdot \nabla (\log n(x,t) - c(x,t))\\  
&=&- \int_{\R^2} n(x,t)|\nabla(\log n(x,t) - c(x,t)+\log(1+|x|^2)) |^2 \ dx - \e\int_{\R^2} |\partial_t c(x,t)|^2\ dx\\  
& &\quad+\int_{\R^2} n(x,t)|\nabla\log(1+|x|^2)) |^2\ dx\ ,
\end{eqnarray*}
where
$$
\int_{\R^2} n(x,t)\left|\nabla\log(1+|x|^2) \right|^2\ dx=\int_{\R^2} \f{4|x|^2}{(1+|x|^2)^2}\ n(x,t)\ dx\le M\ .
$$
As a consequence, $\mathcal E_H(t)$ grows at most linearly in time.\\

Before concluding this preliminary section, we state the minimization lemmas for the entropy and the chemical energy respectively. Moreover, we derive also some bounds for the moments of the density $n$ in term of the entropy production \fer{entropyproduction} that is shown to be locally integrable in time in the sequel (see Theorem \ref{th:aprioriestimate}). 

\begin{lem}[The entropy minimization] 
\label{lem:entropy}
Let $\psi$ be any function such that $e^{ \psi}\in L^1(\R^2)$ and denote $\overline{n}  = M e^{ \psi}\Big(\int_{\R^2} e^{\psi}\ dx\Big)^{-1}$, with $M$ a positive arbitrary constant. Let $E\ :\ L^1_+(\R^2)\to\R\cup\{\infty\}$ be the entropy functional 
$$
E(n;\psi)=\int_{\R^2}\ (n(x)\log n(x)-n(x)\psi(x))\ dx
$$ 
and $RE\ :\ L^1_+(\R^2)\to\R\cup\{\infty\}$
$$
RE(n|\overline{n})=\int_{\R^2}\ n(x)\log\big( n(x)/{\overline n}(x)\big)\ dx
$$ 
the relative (to $\overline n$) entropy. Then, $E(n;\psi)$ and $RE(n|\overline{n})$ are finite or infinite in the same time and for all $n$ in the set ${\cal U}=\{n\in L_+^1(\R^2)\  ,\ \int_{\R^2}n(x)\ dx=M\}$ it holds true that 
\begin{equation}
 E(n;\psi) - E(\overline{n};\psi) = RE(n|\overline{n})\geq 0\ .
 \label{entropyequality}  
\end{equation}
\end{lem}

The entropy minimization Lemma \ref{lem:entropy} is now a classical lemma and the proof can be found for exemple in \cite{CJMTU}, where a more general class of entropy functionals including $E(n;\psi)$ is considered. Anyway, being this lemma of primary importance and  for the sake of completeness, we will give the proof in the appendix \ref{ap:entropy}. 

\begin{lem}[The chemical energy minimization] 
\label{lem:chemicalmin} 
Assume $\alpha\ge0$ and let $f$ be an $L^1(\R^2)$ function such that 
\begin{description}
\item[(A1)] if $\alpha>0$ then $f\ge0$ and $f\log f\in L^1(\R^2)$\ ;
\item[(A2)] if $\alpha=0$ then $f\in L^1(\R^2,\log(1+|x|^2)dx)$, $\int_{\R^2} f(x)\ dx=0$ and $f$ can be decomposed as $f=f_1+f_2$ with $f_1\ge0$, $f_1\in L^1_{loc}(\R^2)$, $f_1\log f_1\in L^1(\R^2)$ and $f_2\in L^2(\R^2)$.
\end{description}
\noindent Finally, let us denote
\begin{equation}
\overline c(x)=\left\{
\begin{array}{ll}
(B_\alpha*f)(x)&\quad\hbox{if }\alpha>0\ ,\\
(E_2*f)(x)&\quad\hbox{if }\alpha=0\ ,
\end{array}
\right.
\label{quasistationarystate}
\end{equation}
where $*$ is the space convolution, $B_\alpha$ denotes the Bessel kernel $
B_\alpha(z) =  \f1{4\pi}\int_0^{+\infty}\f 1t\  e^{-\frac{|z|^2}{4t} - \alpha t}\ dt$ and 
$E_2(z)=-\f1{2\pi}\log|z|$ is the fundamental solution of the Laplace's equation in $\R^2$. Then, $\overline c\in H^1(\R^2)$ when $\alpha>0$ while $|\nabla \overline c|\in L^2(\R^2)$ when $\alpha=0$, with
\begin{equation}
\nabla\overline c(x)=\left\{
\begin{array}{ll}
(\nabla B_\alpha*f)(x)&\quad\hbox{if }\alpha>0\ ,\\
(\nabla E_2*f)(x)&\quad\hbox{if }\alpha=0\ .
\end{array}
\right.
\label{gradc}
\end{equation}
Moreover, if $\alpha>0$ then for all $c\in H^1(\R^2)$, the chemical energy $F_\alpha(c;f)$ is finite and it satisfies 
\begin{equation}
F_\alpha(c;f)- F_\alpha(\overline c;f)=\f12\int_{\R^2}|\nabla (c-\overline c)(x)|^2\ dx+\f\alpha2\int_{\R^2}(c-\overline c)^2(x)\ dx\geq 0\ ,
\label{bo}
\end{equation}
while if $\alpha=0$ then for all $c$ such that $|\nabla c|\in L^2(\R^2)$, the chemical energy $F_0(c;f)$ is finite and it satisfies 
\begin{equation}
F_0(c;f)- F_0(\overline c;f)=\f12\int_{\R^2}|\nabla (c-\overline c)(x)|^2\ dx \geq 0\ .
\label{bo2}
\end{equation}
\end{lem}

Let us observe that whenever $f\in L^2(\R^2)$ in \fer{quasistationarystate} the chemical energy minimization Lemma \ref{lem:chemicalmin} can be obtained easily applying the variational method, at least for $\alpha>0$. However, we want to use here minimal hypotheses on $f$ and therefore the proof becomes a little more technical, expecially when $\alpha=0$. Again, for the sake of completeness, the proof is given in appendix \ref{ap:chemicalminimization}

\begin{lem}[Moment Lemma]
\label{lemma:momentum} 
Let $c$ be a given smooth function. Let $n$ be any non-negative and sufficiently smooth solution of $\partial_tn=\Delta n-\nabla\cdot( n\nabla c)$, with fast decay at infinity and total mass $M$. Then, for any $\delta>0$, we have the following bounds for the evolution of the moments of $n$:
\begin{align}
\int_{\R^2}\!\!n(x,t)\log(1+ |x|^2)\ dx \le & \int_{\R^2}\!\!n_0(x)\log(1+ |x|^2)\ dx+
\f M{2\delta}t \nonumber 
\\ & \qquad  +\f\delta2\!\int_0^t\!\!\int_{\R^2}\! n(x,s)|\nabla(\log n(x,s) - c(x,s)) |^2 \ dxds\ ;
\label{eq:momlog}
\end{align}
\begin{align}
\int_{\R^2} |x|n(x,t)\ dx\le &\int_{\R^2} |x|n_0(x)\ dx+
\f M{2\delta}t \nonumber \\
&\qquad +\f\delta2\int_0^t\!\!\int_{\R^2} n(x,s)|\nabla(\log n(x,s) - c(x,s)) |^2 \ dxds\ ;
\label{eq:mom1}
\end{align}
\begin{equation}
\int_{\R^2} |x|^2n(x,t)\ dx\le 2\int_{\R^2} |x|^2n_0(x)\ dx
+2\ t\int_0^t\!\!\int_{\R^2} n(x,s)|\nabla(\log n(x,s) - c(x,s)) |^2 \ dxds\ .
\label{eq:mom2}
\end{equation}
\end{lem}
\begin{proof} Writing the equation on $n$ as $\partial_t n=\nabla\cdot\big(n\nabla(\log n-c)\big)$, it follows that
\begin{eqnarray*}
 \dfrac{d}{dt}  \int_{\R^2}n(x,t)\phi(x)\ dx & = &  -\int_{\R^2} n(x,t) \nabla \phi(x) \cdot \nabla (\log n(x,t) - c(x,t))\ dx \\
& \leq & \f1{2\delta} \int_{\R^2} \big|\nabla\phi(x)\big|^2 n(x,t)\ dx + \f\delta2\int_{\R^2} n(x,t)|\nabla (\log n(x,t) - c(x,t))|^2\ dx\ ,
\end{eqnarray*}
with $\delta>0$ arbitrary. Taking successively  $\phi(x) = \log(1+|x|^2)$ and $\phi(x) = |x|$, observing that with the first choice of $\phi$ we have
\[ |\nabla\phi(x) | = \left|\frac{2x}{1+|x|^2}\right| \leq 1\ ,\]
and using the mass conservation property, inequalities \fer{eq:momlog} and \fer{eq:mom1} are proved. Inequality \fer{eq:mom2} follows in a similar way since for $\phi(x)=|x|^2$ it holds true that
\begin{eqnarray}
 \dfrac{d}{dt}  \int_{\R^2}|x|^2n(x,t)\ dx & = &  -2\int_{\R^2} n(x,t)\ x \cdot \nabla (\log n(x,t) - c(x,t))\ dx  \nonumber\\ 
&\le& 2\left(\int_{\R^2}|x|^2n(x,t)\ dx\right)^{1/2}
\left(\int_{\R^2} n(x,t)|\nabla (\log n(x,t) - c(x,t))|^2\ dx\right)^{1/2}\ .
\nonumber\\
\label{estmom2}
\end{eqnarray}
Thereby, integrating \fer{estmom2}, we obtain
\begin{eqnarray*}
\int_{\R^2}|x|^2n(x,t)\ dx&\le&\left[\left(\int_{\R^2}|x|^2n_0(x)\ dx\right)^{1/2}
+\int_0^t\left(\int_{\R^2} n(x,s)|\nabla (\log n(x,s) - c(x,s))|^2\ dx\right)^{1/2}ds\right]^2\\
&\le& 2\int_{\R^2} |x|^2n_0(x)\ dx
+2\ t\int_0^t\int_{\R^2} n(x,s)|\nabla(\log n(x,s) - c(x,s)) |^2 \ dx\ ds\ ,
\end{eqnarray*}
and the lemma is proved.

\end{proof}

The proof of the previous lemma is based uniquely on the equation on $n$ and on the specific expression of the weight function defining the moment. The evolution followed by $c$ doesn't play any role and the lemma holds true also for the (PKS) system with non-negative coefficients $\e$, $\eta$ and $\alpha$.  Of course, one can estimate the evolution of other moments than those considered in the lemma. Here we have considered the most useful and used. In particular, the local in time bound of the weighted $L^1(\R^2,\log(1+ |x|^2)dx)$ norm of $n$ together with the following lemma will be of primary importance to obtain the key {\em a priori} bounds and the key equi-integrability of $n$ giving the global existence with a regularizing effect.

\begin{lem}
\label{lem:controlmass}
Let $\psi$ be any function such that $e^{\psi}\in L^1(\R^2)$ and let $f$ be a non-negative function such that $(f\ \bone_{\{f\le1\}})\in L^1(\R^2)\cap  L^1(\R^2,|\psi(x)|dx)$. Then, there exists a constant $C$ such that
\begin{equation}
\int_{\R^2} f(x) (\log f(x))_-\ dx \leq C-\int_{\{f\le1\}} f(x)\psi(x)\ dx\ .
\label{eq:masscontrol}
\end{equation}
\end{lem}
\begin{proof} 
Let us denote $v=f\ \bone _{\{f\le1\}}$, $m=\int_{\R^2}v(x)\ dx<\infty$ and $\overline v(x)=me^{\psi(x)}(\int_{\R^2}e^{\psi(x)}\ dx)^{-1}$. Then, by the entropy minimization Lemma \ref{lem:entropy} we have
\begin{eqnarray*}
\int_{\R^2}(v(x)\log v(x)-v(x)\psi(x))\ dx&\ge&
\int_{\R^2}(\overline v(x)\log\overline v(x)-\overline v(x)\psi(x))\ dx\\
&=&m\log m-m\log\left(\int_{\R^2}e^{\psi(x)}\ dx\right)=C\ ,
\end{eqnarray*}
which gives \eqref{eq:masscontrol} thanks to the identity $\int_{\R^2} f(x) (\log f(x))_-\ dx =-\int_{\R^2} v(x)\log v(x)\ dx.$

\end{proof}

\section{\emph{A priori} estimates from the Moser-Trudinger-Onofri inequality}
\label{sec:Onofri}
In \cite{Ono} Onofri obtained the following sharp inequality on the sphere $\Sph^2$
\begin{equation}
\int_{\Sph^2} e^{v(s)}\ ds \le  \exp\left\{\int_{\Sph^2}\left(v(s)+\f14\vert\nabla_0v(s)\vert^2\right)\ ds\right\} \ ,
\label{eq:OnoS2}
\end{equation}
for all functions $v\in L^1(\Sph^2,ds)$ such that $\vert\nabla_0v\vert\in L^2(\Sph^2,ds)$. Here, $ds$ is the uniform normalized surface measure on $\Sph^2$ so that $\int_{\Sph^2}\ ds=1$. Moreover $\nabla_0$ is the covariant gradient with respect to the metric $ds_0^2=d\theta^2+\sin^2\theta\ d\phi^2$, $(\theta,\phi)$ being the polar coordinates, i.e.
$$
|\nabla_0 v|^2=\left(\f{\partial v}{\partial\theta}\right)^2+(\sin\theta)^{-2}\left(\f{\partial v}{\partial\phi}\right)^2\ .
$$
With the help of the stereographic projection ${\cal S}$, the same inequality can be stated equivalently on $\R^2$ as follows
\begin{lem}[Onofri's inequality in $\R^2$]
\label{lem:Onofri}
Let $H$ be defined as in \fer{H}. Then 
\begin{equation}
\int_{\R^2}e^{u(x)}H(x)\ dx\leq\exp\left\{\int_{\R^2}u(x)H(x)\ dx+\f1{16\pi}\int_{\R^2}\vert\nabla u(x)\vert^2\ dx\right\}\ ,
\label{eq:OnoR2}
\end{equation}
for all functions $u\in L^1(\R^2,H(x)dx)$ such that $\vert\nabla u(x)\vert\in L^2(\R^2,dx)$.
\end{lem}
\begin{proof} 
It is sufficient to apply the Onofri's inequality \fer{eq:OnoS2} to the function $e^{\tilde u}$ with $\tilde u=u\circ {\cal S}^{-1}$ and we get
\begin{eqnarray*}
\int_{\R^2}e^{u(x)}H(x)\ dx=\int_{\Sph^2} e^{\tilde u(s)}\ ds &\le& 
\exp\left\{\int_{\Sph^2}\Big( \tilde u(s)+\f14\vert\nabla_0\tilde u(s)\vert^2\Big)\ ds\right\}\\ 
&=&\exp\left\{ \int_{\R^2}u(x)H(x)\ dx+\f1{16\pi}\int_{\R^2}\vert\nabla u(x)\vert^2\ dx\right\}\ ,
\end{eqnarray*}
since $\vert\nabla_0\tilde u(s)\vert^2=(4\pi H(x))^{-1} \vert\nabla u(x)\vert^2$ and $ds= H(x)dx$.

\end{proof}

Thanks to inequality \fer{eq:OnoR2}, we are now able to follow the first strategy giving {\em a priori} estimates, namely the minimization of $\mathcal E(t)$ with respect to $n$. More precisely, first we apply the entropy minimization  Lemma \ref{lem:entropy} with $\psi=(1+\delta)c+\log H$ and then we  make use of the Onofri's inequality \fer{eq:OnoR2} with $u=(1+\delta)c$. Working in this way, we are able to obtain the optimal threshold value, namely $8\pi$, for the mass $M$. Observe that this is the same threshold as the one obtained in \cite{BDP} for the parabolic-elliptic (PKS) system with $\e=\alpha=0$ over $\R^2$, as it would be expected. Moreover, we optimize the result given in \cite{N01}, where the author obtained the global existence of non-negative solutions $(n,c)$ of \fer{KS} over $\R^2$ under the smallness condition $M<4\pi$. This is due to the fact that in \cite{N01} the author use a Brezis-Merle type inequality for the heat equation on $\R^2$ instead of a Moser-Trudinger type inequality as \fer{eq:OnoR2}.

We give the announced estimates first formally in the following Theorem. The procedure leading to the rigorous existence result will be given later in Section \ref{Sec:Globalexistence}.
\begin{thm}[{\em A priori} estimates]
Assume $\e>0$ and $\alpha\ge0$. Under the same hypotheses (H1)-(H2)-(H3) as in Theorem \ref{the:existence} and for subcritical mass $M<8\pi$, let $(n,c)$ be a non-negative solution of \fer{KS}, supposed sufficiently smooth with fast decay as $|x|\to+\infty$. Then, $n\log n\in L^\infty(0,T;L^1(\R^2))$, $n\in L^\infty(0,T;L^1(\R^2,\log(1+|x|^2)dx))$, for any $T>0$,  and the  following {\it a priori} estimates hold true for all $t>0$ :
\begin{enumerate}[(i)]
\item 
$ \int_{\R^2} n(x,t)(\log n(x,t))_+\ dx \le C(1+t)\ ;$
\item 
$\int_{\R^2} n(x,t)c(x,t)\ dx\le C(1+t)\ ;$
\item 
$\| c(t)\|^2_{H^1(\R^2)}\le C(1+t)\quad$ if $\alpha>0\ $, \  and
$\| \nabla c(t)\|^2_{L^2(\R^2)}\le C(1+t)\quad$ if $\alpha=0\ ;$
\item
$\mathcal E(t)\ge -C(1+t)\ ;$
\item 
$\mathcal E_H(t)\ge -C(M,\alpha)\quad$ if $\alpha>0\ $,\ and 
$\quad\mathcal E_H(t)\ge -C(1+t)\quad$ if $\alpha=0\ ;$
\item
$\int_0^t\int_{\R^2} n(x,s) |\nabla (\log n(x,s) - c(x,s)) |^2\ dxds\le\mathcal E(0)+
2\int_{\R^2}n_0(x) \log(1+|x|^2)dx+C(1+t) \ ;$
\item 
$\e \int_0^t\int_{\R^2} |\partial_tc(x,s)|^2\ ds\le C(1+t)\ .$
\end{enumerate}
\label{th:aprioriestimate}
\end{thm}
\begin{proof} 
Let us consider first the case $\alpha>0$. Let $\delta>0$ and $\tilde\delta>0$ to be chosen later and let us define $\displaystyle\overline{n}~=~M e^{(1+\delta) c(x,t)}H(x)\left(\int_{\R^2}e^{(1+\delta) c(x,t)}H(x)\ dx\right)^{-1}$. We observe that, thanks to the Onofri's inequality \fer{eq:OnoR2}, $c(t)\in H^1(\R^2)$ is sufficient in order to have $\overline n$ well defined and $\overline n(t)\in L^1(\R^2)$. Then, we can apply the Entropy Lemma \ref{lem:entropy} with $\psi=(1+\delta)c+\log H$ to obtain
\begin{eqnarray}
 E(n;(1+\delta) c+\log H) &\ge&  E(\overline{n};(1+\delta) c+\log H) \nonumber 
\\ &=& M\log M-M\log\left(\int_{\R^2}e^{(1+\delta) c(x,t)}H(x)\ dx\right).
 \label{mini}
\end{eqnarray}
Furthermore, applying Lemma \ref{lem:Onofri} with $u = (1+\delta )c$ to the last term in the right hand side of \fer{mini}, we have 
for the modified free energy functional $\mathcal E_H(t)$ \fer{eq:feH},
\begin{eqnarray}
\mathcal E_H(t)&=&  E(n;(1+\delta)c+\log H)+\delta\int_{\R^2} n(x,t)c(x,t)\ dx+ \f12\int_{\R^2} |\nabla c(x,t)|^2\ dx + \f\alpha2 \int_{\R^2} c^2(x,t)\ dx
\nonumber\\ 
&\ge&M\log M +\f12\left(1-M\f{(1+\delta)^2}{8\pi}\right)\int_{\R^2}\vert\nabla c(x,t)\vert^2\ dx-M(1+\delta)\int_{\R^2}c(x,t)H(x)\ dx
\nonumber\\  
&\quad&\quad\quad+\delta\int_{\R^2} n(x,t)c(x,t)\ dx+ \f\alpha2\int_{\R^2} c^2(x,t)\ dx
\label{eq:cH}\\ 
&\ge& \f12\left(1-M\f{(1+\delta)^2}{8\pi}\right)\int_{\R^2}\vert\nabla c(x,t)\vert^2\ dx+\left(\f\alpha2-M(1+\delta)\f{\tilde\delta}2\right)\int_{\R^2} c^2(x,t)\ dx +M\log M
\nonumber\\  
&\quad&\quad\quad+\delta\int_{\R^2} n(x,t)c(x,t)\ dx-M(1+\delta)\f1{2\tilde\delta}\int_{\R^2} H^2(x)\ dx\ .
\label{estEH}
\end{eqnarray}
Next, we choose $\delta>0$ small enough such that $M<\dfrac{8\pi}{(1+\delta)^2}$ and $\tilde\delta>0$ such that $\alpha>M(1+\delta)\tilde\delta$. This is possible because $M$ is less than the critical mass value $8\pi$. As a consequence, since $\mathcal E_H(t)$ grows at most linearly, \fer{estEH} gives us (ii), (iii) and (v). 

When $\alpha=0$, we have to estimate differently the term $\int_{\R^2} c(x,t)H(x)\ dx$ in \eqref{eq:cH}. That can be done using the following identity,
\begin{equation}
\int_{\R^2}c(x,t) \ dx=\dfrac1\e {Mt} + \int_{\R^2}c_0(x)\ dx
\label{eq:L1normc}
\end{equation}
and the fact that $H$ is bounded. Again, whenever $|\nabla c(\cdot,t)|\in L^2(\R^2)$, we have that $\overline n$ is well defined, $\overline n(\cdot,t)\in L^1(\R^2)$ and we can obtain (ii), (iii) and (v) as before.

From now on, let $\alpha\ge0$. Proposition \ref{prop:fe} , the Moment Lemma \ref{lemma:momentum} and the definition \fer{H} of $H$ give us the following estimate 
\begin{align*}
&\int_0^t\int_{\R^2} n(x,s) |\nabla (\log n(x,s)-c(x,s)) |^2\ dx ds
\\ &\qquad \le\mathcal E(0)-\mathcal E(t)=\mathcal E(0)-\mathcal E_H(t)-\int_{\R^2}n(x,t)\log H(x)\ dx  \\  
&\qquad =\mathcal E(0)-\mathcal E_H(t)+M\log\pi+2\int_{\R^2}n(x,t)\log(1+|x|^2)\ dx  \\  
&\qquad \le\mathcal E(0)-\mathcal E_H(t)+M\log\pi+2\int_{\R^2}n_0(x)\log(1+|x|^2)\ dx
+\f1\delta Mt  \\  
&\qquad \qquad +\delta\int_0^t\int_{\R^2} n(x,s) |\nabla (\log n(x,s) - c(x,s)) |^2\ dx ds\ . 
\end{align*}
Choosing $\delta<1$ and using the lower bound (v) on $\mathcal E_H(t)$, we obtain (vi). As a consequence of (vi) and of the Moment Lemma \ref{lemma:momentum} again, we have that $n\in L^\infty(0,T;L^1(\R^2,\log(1+|x|^2)dx))$. 

In the same way, we obtain (vii) and the lower bound (iv) for the free energy $\mathcal E(t)$ since $\e \int_0^t\int_{\R^2} |\partial_tc(s)|^2\ ds\le\mathcal E(0)-\mathcal E(t)$ and $\mathcal E(t)=\mathcal E_H(t)+\int_{\R^2}n(x,t)\log H(x)\ dx$. 

To conclude it remains to prove (i) and that $n\log n\in L^\infty(0,T;L^1(\R^2))$. This is a straightforward consequence of Lemma \ref{lem:controlmass}. Indeed, from the free energy definition \fer{eq:fe} and the previous estimates it follows that
$$
-C(1+t)\le\int_{\R^2} n(x,t)\log n(x,t)\ dx\le C(1+t)\ .
$$
Next, appling Lemma \ref{lem:controlmass} to $n$ with $\psi(x)=-(1+\delta)\log(1+|x|^2)$ and arbitrary $\delta>0$ in order to have $e^{\psi}\in L^1(\R^2)$, we have
\begin{eqnarray*}
\int_{\R^2} n(x,t)(\log n(x,t))_+ \ dx&=&\int_{\R^2} n(x,t)\log n(x,t)\ dx + \int_{\R^2} n(x,t)(\log n(x,t))_- \ dx\\
&\le& C(1+t) + (1+\delta) \int_{\R^2} n(x,t)\log (1+|x|^2)\ dx + C\le C(1+t) 
\end{eqnarray*}
and we have obtained (i). Finally, the identity
$$
\int_{\R^2}| n\log n|\ dx=\int_{\R^2} n\log n\ dx + 2 \int_{\R^2} n(x,t)(\log n(x,t))_- \ dx \ ,
$$
gives us that $n\log n\in L^\infty(0,T;L^1(\R^2))$.
\end{proof}

\begin{rem} [The critical case $M=8\pi$]. We cannot afford the critical mass value $M=8\pi$ in the previous Theorem. Indeed this necessarily leads to $\delta=0$ in \fer{estEH}. Then one can prove successively estimates (v), $\|c(t)\|^2_{L^2(\R^2)}\le C(1+t)$ if $\alpha>0$, (vi), $n\in L^\infty(0,T;L^1(\R^2,\log(1+|x|^2)dx))$, (iv), (vii), but we can't obtain the fundamental estimates (i), (ii) and (iii).
\label{rem:8pi}
\end{rem}

\begin{rem}[The case $\alpha=0$].
\label{rem:weakstrongc0} When $\alpha=0$, in Theorem \ref{th:aprioriestimate} we can use the necessary and sufficient hypothesis on $c_0$ to apply the Onofri's inequality \fer{eq:OnoR2}, i.e. $c_0\in L^1(\R^2,H(x)dx)$ instead of $c_0\in L^1(\R^2)$ and $|\nabla c_0|\in L^2(\R^2)$. However, we can obtain only estimates exponentially increasing in time. Indeed, using the identity
$$
\e\f d{dt}\int_{\R^2}c(x,t)H(x)\ dx=\int_{\R^2}c(x,t)\Delta H(x)\ dx+\int_{\R^2}n(x,t)H(x)\ dx
$$
and the computation $\Delta H(x)=\f8\pi\left[\f{2|x|^2-1}{(1+|x|^2)^4}\right]$, we easily have
$$
\e\f d{dt}\int_{\R^2}c(x,t)H(x)\ dx\le4\int_{\R^2}c(x,t) H(x)\ dx+\frac M\pi\ .
$$
Then, 
\begin{equation}
\int_{\R^2}c(x,t)H(x)\ dx\le e^{\f4\e t}\left(\int_{\R^2}c_0(x)H(x)\ dx+\f M{4\pi}\right)\ .
\label{estcH}
\end{equation}
Injecting \fer{estcH} in \fer{eq:cH} instead of \fer{eq:L1normc}, we have the assertion.

\noindent On the other hand, under the stronger hypothesis $c_0\in H^1(\R^2)$, Theorem \ref{th:aprioriestimate} holds true in the case $\alpha=0$ with the additional estimate $\int_{\R^2}c^2(x,t)\ dx\le C(1+t+t^2)$ that follows by
$$
\f\e2\int_{\R^2}c^2(x,t)\ dx\le\f\e2\int_{\R^2}c_0^2(x)\ dx+\int_0^t\int_{\R^2}n(x,t)c(x,t)\ dx\ .
$$
\end{rem}

\begin{rem}[The case $\alpha=0$, corrected solutions].
\label{rem:correctedchemical}
In the peculiar case $\alpha=0$, one can develop another sort of solutions based on the following chemical deviation 
$$
u(x,t) := c(x,t) - \f M{8\pi}\log H(x),
$$
with $M<8\pi$. With this change of variable, it is sufficient to assume initially $c_0\in L^1(\R^2, H(x)dx)$ and $\nabla \left(c_0 + \frac M{4\pi}\log (1+|x|^2)\right)\in L^2(\R^2)$. The second hypothesis is hardly biologically relevant because it breaks the non-negativity assumption on $c_0$. However, it allows us to solve the technical difficulty arising from evaluating the integral $\int_{\R^2}c(x,t) H(x)\ dx$ in \eqref{eq:cH} and to obtain the same estimates as in Theorem \ref{th:aprioriestimate} (except for (iii) where $\nabla u$ has to be read instead for $\nabla c$) in the following way. 

The function $\log H$ satisfies the remarkable equation
\begin{equation} 
-\Delta \log H(x)=2\Delta \log (1+|x|^2)=8\pi H(x) 
\label{eq:cancellation} \ .
\end{equation}
Therefore, the pair $(n,u)$ verifies the system
\begin{equation}
\left\{\begin{array}{rcl}
\partial_t n& =& \Delta  n - \nabla\cdot( n\nabla( u +\f M{8\pi}\log H )) \\
\e\partial_t u&= &\Delta u +  n - M H \ ,
\end{array}
\right. 
\label{eq:KScorr}
\end{equation}
which admits the following free energy
\begin{equation}
\begin{split} 
\tilde{\mathcal E}(t) = \int_{\R^2}n(x,t) \log n(x,t)\ dx - \int_{\R^2}n(x,t) u(x,t)\ dx + \frac12 \int_{\R^2}|\nabla u(x,t)|^2\ dx \\  
- \f M{8\pi}\int_{\R^2}n(x,t) \log H(x) \ dx + M \int_{\R^2}u(x,t) H(x)\ dx\ .
\end{split}
\label{eq:energycorr} 
\end{equation}
It is worth noticing that the free energy $\tilde{\mathcal E}(t)$, up to an additional constant, is exactly the original free energy $\mathcal E(t)$ corresponding to the variables $(n,c)$ and that
\begin{equation}
\dfrac{d}{dt}\tilde{\mathcal E}(t) = - \int_{\R^2} n(x,t)|\nabla(\log n(x,t) - c(x,t)) |^2 \ dx - \e\int_{\R^2} |\partial_t c(x,t)|^2\ dx\ .
\label{eq:evEcorr}
\end{equation}
Then, reasoning as before, we write
\begin{align} 
 \tilde{\mathcal E}(t) =&\delta \int_{\R^2}n(x,t) \log n(x,t)dx +(1-\delta)E\left(n;\f u{(1-\delta)}+\log H\right) \nonumber \\ &\qquad + \left(1-\delta-\f M{8\pi}\right) \int_{\R^2}n(x,t) \log H(x)dx
+\f12\int_{\R^2}|\nabla u(x,t)|^2\ dx+M \int_{\R^2}u(x,t) H(x)\ dx\ ,\quad\quad
\label{eq:energycorr2} 
\end{align}
with some degree of freedom $0<\delta<1$ to be chosen later. We next apply the Entropy Lemma \ref{lem:entropy} and the Onofri's inequality \eqref{eq:OnoR2} to obtain
\begin{eqnarray*}
E\left(n;\f u{(1-\delta)}+\log H\right)&\ge& M\log M
-M\log\left(\int_{\R^2}e^{u(x,t)/(1-\delta)}H(x)dx\right)\\
&\ge&M\log M-\f M{(1-\delta)}\int_{\R^2}u(x,t)H(x)\ dx-\f M{16\pi(1-\delta)^2}\int_{\R^2}|\nabla u(x,t)|^2\ dx\ .
\end{eqnarray*}
Plugging this estimate into the corrected free energy \eqref{eq:energycorr2}, we deduce
\begin{eqnarray*}
\tilde{\mathcal E}(t) & \geq & \delta \int_{\R^2}n(x,t) \log n(x,t)\ dx + \left(1-\delta-\f M{8\pi}\right) \int_{\R^2}n(x,t) \log H(x)\ dx \\
&&\qquad + \frac12\left( 1 - \dfrac{M}{8\pi(1-\delta)} \right) \int_{\R^2}|\nabla u(x,t)|^2\ dx +(1-\delta)M\log M\ .
\end{eqnarray*}
Choosing $0<\delta<1$ small enough to ensure $\left(1 - \dfrac{M}{8\pi(1-\delta)}\right) >0$, the bootstrap argument goes as previously. From the definition of $H$, the Moment Lemma \ref{lemma:momentum} and the evolution equation \fer{eq:evEcorr} of $\tilde{\mathcal E}(t)$, we obtain
\begin{align*} 
&\delta \int_{\R^2}n(x,t) \log n(x,t)\ dx + \frac12\left( 1 - \dfrac{M}{8\pi(1-\delta)} \right) \int_{\R^2}|\nabla u(x,t)|^2\ dx \\
&\qquad \qquad\qquad\qquad\qquad+ (\f M{8\pi}+\delta) \int_0^t \int_{\R^2} n(x,s) |\nabla (\log n(x,s)-c(x,s)) |^2\ dx ds \\
&\qquad \qquad\qquad \leq  \tilde{\mathcal E}(0) + 2(1-\delta-\f M{8\pi}) \int_{\R^2}n_0(x) \log(1+|x|^2)\ dx +  (1-\delta-\f M{8\pi})M t \ .
\end{align*}
Finally, it remains to apply Lemma \ref{lem:controlmass} as done in Theorem \ref{th:aprioriestimate} to control the contribution $\int_{\R^2}n(\log n)_-\ dx$. 
\end{rem}

\section{\emph{A priori} estimates from the logarithmic HLS inequality}
\label{sec:log HLS}
In this section we follow the strategy dual to the previous one, starting from the  minimization of the modified free energy functional $\mathcal E_H(t)$  with respect to $c$ and giving the same a priori estimates. The role of the Onofri's inequality \fer{eq:OnoR2} will be played by the logarithmic HLS inequality below when $\alpha=0$ and by its generalization to the Bessel kernel $B_\alpha$ when $\alpha>0$ (see Lemma \ref{lem:Besselinequality}). Therefore, the cases $\alpha>0$ and $\alpha=0$ have to be treated separately again.

\begin{lem}[Logarithmic Hardy-Littlewood-Sobolev inequality in $\R^2$]
For all non-negative functions $f\in L^1(\R^2)$ such that $\int_{\R^2} f(x)\ dx = M$, $\int_{\R^2} f(x) \log(1+|x|^2)\ dx<\infty$ and $\int_{\R^2}f(x)\log f(x)\ dx$ is finite, the following inequality holds true
\begin{equation}
-\int_{\R^2}\int_{\R^2} f(x)\log|x-y|f(y)\, dx dy \leq
 \f M2\int_{\R^2}f(x)\log f(x)\,dx+C(M)\ .
\label{eq:logHLS}
\end{equation}
\label{lem:logHLS}
\end{lem}

The proof of Lemma \ref{lem:logHLS} can be found for example in \cite{CarlenLoss}. Let us observe here that the logarithmic HLS inequality \fer{eq:logHLS} written equivalently on $\Sph^2$ is dual to the Onofri's inequality \fer{eq:OnoS2} in the sense that extremal functions for one inequality determine extremal functions for the other inequality (see \cite{Beckner}, \cite{CarlenLoss}). As an exemple of this kind of duality results, we show in the appendix \ref{duality} how one can obtain directly inequality \fer{eq:Besselinequality} from the Onofri's inequality on $\R^2$ \fer{eq:OnoR2}.

The generalization of inequality \fer{eq:logHLS} to the Bessel kernel $B_\alpha$ (whose definition is given in Lemma \ref{lem:chemicalmin}) is the following Lemma \ref{lem:Besselinequality}. One possible interpretation why an additional ``logarithmic" momentum appears in \eqref{eq:Besselinequality} is that the homogeneity of the logarithmic kernel $-\log|z|$ has been broken in the Bessel potential $B_\alpha$ (consider for instance the dilation $f_\lambda(x) = \lambda f(\lambda x)$). It shares indeed  the same singularity for small $|z|$ but  for large $|z|$ $B_\alpha$ is a positive exponentially decreasing function. For other extensions of Lemma~\ref{lem:logHLS} see \cite{Beckner}.

\begin{lem}
For all non-negative functions $f\in L^1(\R^2)$ such that $\int_{\R^2} f(x)dx = M$, $\int_{\R^2}f(x)\log f(x)dx$ is finite and $\int_{\R^2} f(x) \log(1 +|x|^2)\ dx<\infty$, the following inequality holds true
\begin{equation}
\begin{split}
\int_{\R^2}\int_{\R^2} f(x)B_\alpha(x-y)f(y)\, dx\, dy \leq
 \f M{4\pi}\int_{\R^2}f(x)\log f(x)\,dx \\ + \f M{2\pi}\int_{\R^2}f(x)\log(1+|x|^2)\,dx+C(M).
\label{eq:Besselinequality}
\end{split}
\end{equation}
\label{lem:Besselinequality}
\end{lem}

\begin{proof} The Bessel kernel $B_\alpha$ is a positive radial decreasing function such that $B_\alpha(z)\to+\infty$ when $|z|\to0$ as $\frac1{2\pi}\log\left(\f1{|z|}\right) + O(1)$. We only prove below the upper side of this asymptotic estimate, which is sufficient for our purpose. Indeed,  following \cite{Dono}, $B_\alpha$ can be written as the sum of the three integrals $\int_0^{r^2}+\int_{r^2}^1+\int_1^{+\infty}\ $, with $r=|z|<1$.  For the first and third integrals we have respectively
$$
\f1{4\pi}\int_0^{r^2}\f 1t\  e^{-\frac{r^2}{4t} - \alpha t}\ dt\le 
\f1{4\pi}\int_0^{r^2}\f 1t\  e^{-\frac{r^2}{4t}}\ dt=\f1{4\pi}\int_{\f14}^{+\infty}\f{e^{-y}}y\ dy<\infty
$$
and 
$$
\f1{4\pi}\int_1^{+\infty}\f 1t\  e^{-\frac{r^2}{4t} - \alpha t}\ dt\le 
\f1{4\pi}\int_1^{+\infty}e^{- \alpha t}\ dt<\infty\ .
$$
The second integral satisfies
\begin{equation}
\f1{4\pi} e^{-(\alpha+\f14)}\log\left(\f1{r^2}\right)\le\f1{4\pi}\int_{r^2}^1\f 1t\  e^{-\frac{r^2}{4t} - \alpha t}\ dt\le
\f1{4\pi}e^{-r^2(\alpha+\f14)}\log\left(\f1{r^2}\right)\ .
\label{est}
\end{equation}

As a consequence of \fer{est}, $\lim_{|z|\to0}B_\alpha(z)=+\infty$ and
\begin{eqnarray*}
\int_{\R^2}\int_{\R^2} f(x)B_\alpha(x-y)f(y)\, dx\, dy &\leq&
\int_{\R^2}\int_{|x-y|\le1}f(x)\left(C-\f1{2\pi}\log|x-y|\right)f(y)\, dx\, dy\\
&\ &\quad\quad\quad+\int_{\R^2}\int_{|x-y|>1}f(x)B_\alpha(x-y)f(y)\, dx\, dy \\
&\le&CM^2-\f1{2\pi}\int_{\R^2}\int_{\R^2} f(x)\log|x-y|f(y)\, dx\, dy \\
&\ &\quad\quad\quad+\f1{2\pi}\int_{\R^2}\int_{|x-y|>1}f(x)\log|x-y|f(y)\, dx\, dy\ .
\end{eqnarray*}
Using the inequality 
\begin{equation}
\log|x-y|\le\f12\log2+\f12\log(1+|x|^2)+\f12\log(1+|y|^2)\ ,
\label{loginequality}
\end{equation}
it follows that
$$
\f1{2\pi}\int_{\R^2}\int_{|x-y|>1}f(x)\log|x-y|f(y)\, dx\, dy\le CM^2+\f M{2\pi}\int_{\R^2}f(x)\log(1+|x|^2)\ dx\ .
$$
Applying the logarithmic HLS inequality \fer{eq:logHLS}, we obtain \fer{eq:Besselinequality}.

\end{proof}

\begin{proof}[Proof of Theorem \ref{th:aprioriestimate}] We first consider the case $\alpha>0$. Let $0<\delta<1$ to be chosen later and let $\overline c$ be the quasi-stationary state \fer{quasistationarystate} corresponding to $(1-\delta)^{-1}n$, i.e. $\overline c(x,t)=\f1{(1-\delta)}(B_\alpha*n(t))(x)$. Then, by the minimization Lemma \ref{lem:chemicalmin}, the chemical energy \fer{chemicalenergy} satisfies
\begin{equation}
F_\alpha\left(c;\f n{1-\delta}\right)\ge F_\alpha\left(\overline c;\f n{1-\delta}\right) 
=-\f1{2(1-\delta)}\int_{\R^2}n(x,t)\overline c(x,t)\ dx\ ,\quad\alpha>0\ .
\label{minic}
\end{equation}
This is possible as soon as $n(t)\in L^1(\R^2)$ and $n(t)\log n(t)\in L^1(\R^2)$. 
Next, from \fer{minic} and Lemma \ref{lem:Besselinequality} we obtain
\begin{eqnarray*}
\mathcal E(t)&=&\int_{\R^2}n(x,t)\log n(x,t)\ dx+(1-\delta) F_\alpha\left(c;\f n{1-\delta}\right)+\delta\left(\f12\int_{\R^2}|\nabla c(x,t)|^2\ dx+ \f\alpha2 \int_{\R^2} c^2(x,t)\ dx\right)\\
&\ge&\int_{\R^2}n(x,t)\log n(x,t)\ dx
-\f1{2(1-\delta)}\iint_{\R^2\times\R^2}n(x,t)B_\alpha(x-y)n(y,t)\ dx\ dy\\
&\ &\quad\quad\quad\quad\quad\quad\quad\quad\quad\quad\quad\quad\quad\quad
+\delta\left(\f12\int_{\R^2}|\nabla c(x,t)|^2\ dx+ \f\alpha2 \int_{\R^2} c^2(x,t)\ dx\right)\\
&\ge&\left(1-\f M{8\pi(1-\delta)}\right)\int_{\R^2}n(x,t)\log n(x,t)\ dx
-\f M{4\pi(1-\delta)}\int_{\R^2}n(x,t)\log (1+|x|^2)\ dx-\f 1{1-\delta} C(M)\\
&\ &\quad\quad\quad\quad\quad\quad\quad\quad\quad\quad\quad\quad\quad\quad
+\delta\left(\f12\int_{\R^2}|\nabla c(x,t)|^2\ dx+ \f\alpha2 \int_{\R^2} c^2(x,t)\ dx\right)\ .
\end{eqnarray*}
Moreover, using the definition \fer{H} of $H$ it follows that
\begin{eqnarray*}
\mathcal E(t)&\ge&\left(1-\f M{8\pi(1-\delta)}\right)\int_{\R^2}n(x,t)\log n(x,t)\ dx
+\f M{8\pi(1-\delta)}\int_{\R^2}n(x,t)\log H(x)\ dx
\\
&\ &\quad\quad\quad\quad\quad\quad\quad\quad\quad\quad\quad\quad
+\delta\left(\f12\int_{\R^2}|\nabla c(x,t)|^2\ dx+ \f\alpha2 \int_{\R^2} c^2(x,t)\ dx\right)-\f 1{1-\delta} C(M)\ .
\end{eqnarray*}
Therefore, the modified free energy $\mathcal E_H(t)$ verifies
\begin{eqnarray*}
\mathcal E_H(t)&\ge&\left(1-\f M{8\pi(1-\delta)}\right)\int_{\R^2}n(x,t)\log\left(\f {n(x,t)}{H(x)}\right)\ dx \\
&\ &\quad\quad\quad\quad\quad\quad\quad\quad\quad\quad\quad\quad
+\delta\left(\f12\int_{\R^2}|\nabla c(x,t)|^2\ dx+ \f\alpha2 \int_{\R^2} c^2(x,t)\ dx\right)-\f 1{1-\delta} C(M)
\end{eqnarray*}
and the Theorem follows from the above estimate. Indeed, since
$$
\int_{\R^2}n(x,t)\log\left(\f {n(x,t)}{H(x)}\right)\ dx=\int_{\R^2}\left(\f {n(x,t)}{H(x)}\right)\log\left(\f {n(x,t)}{H(x)}\right)H(x)\ dx\ge-e^{-1}\ ,
$$
choosing $0<\delta<1$ such that $M<8\pi(1-\delta)$ and using the facts that $\mathcal E_H(t)$ grows at most linearly, estimates (iii) and (v) with $\alpha>0$ of Theorem \ref{th:aprioriestimate} follow. The remaining estimates of Theorem \ref{th:aprioriestimate}  follow as in section  \ref{sec:Onofri}.

To conclude, let us consider the case $\alpha=0$. In this case, the minimization principle \fer{minic} does not hold true since $\nabla \overline{c}(t) = \nabla E_2*n(t)$ does not lie in $L^2(\R^2)$ in general. Nevertheless we can obtain the required estimates by considering the corrected quasi-stationary state $\overline c(x,t)=\f1{(1-\delta)}(E_2*(n(t)-MH))(x)$, whose gradient does belong to $L^2(\R^2)$ by Lemma \ref{lem:chemicalmin}. Indeed, choosing $f_1=n$ and $f_2=-MH$, the hypothesis of Lemma \ref{lem:chemicalmin} are satisfied as soon as $n(t)\in L^1(\R^2)$, $n(t)\log n(t)\in L^1(\R^2)$ and $n(t)\in L^1(\R^2,\log(1+|x|^2)dx)$. Hence,
\begin{eqnarray*}
F_0\left(c;\f 1{1-\delta}(n-MH)\right)&\ge& F_0\left(\overline c;\f 1{1-\delta}(n-MH)\right) \\
&=&-\f1{2(1-\delta)}\int_{\R^2}(n(x,t)-MH(x))\overline c(x,t)\ dx\ .
\end{eqnarray*}
Then, acting as before, we obtain 
\begin{eqnarray*}
\mathcal E(t)&=&\int_{\R^2}n(x,t)\log n(x,t)\ dx+(1-\delta) F_0\left(c;(1-\delta)^{-1}(n-MH)\right)+\f\delta2\int_{\R^2}|\nabla c(x,t)|^2\ dx\\
&\ &\quad\quad\quad\quad\quad\quad\quad\quad\quad\quad\quad\quad-M\int_{\R^2}c(x,t)H(x)\ dx\\
&\ge&\int_{\R^2}n(x,t)\log n(x,t)\ dx+\f1{4\pi(1-\delta)}\iint_{\R^2\times\R^2}(n(x,t)-MH(x))\log|x-y|(n(y,t)-MH(y))\ dx dy\\
&\ &\quad\quad\quad\quad\quad\quad\quad\quad\quad\quad\quad\quad
+\f\delta2\int_{\R^2}|\nabla c(x,t)|^2\ dx-M\int_{\R^2}c(x,t)H(x)\ dx\\
&=&\int_{\R^2}n(x,t)\log n(x,t)\ dx+\f1{4\pi(1-\delta)}\iint_{\R^2\times\R^2}n(x,t)\log|x-y|n(y,t)\ dx dy\\
&\ &\quad\quad-\f M{2\pi(1-\delta)}\iint_{\R^2\times\R^2}H(x)\log|x-y|n(y,t)\ dx dy
+\f {M^2}{4\pi(1-\delta)}\iint_{\R^2\times\R^2}H(x)\log|x-y|H(y)\ dx dy\\
&\ &\quad\quad
+\f\delta2 \int_{\R^2}|\nabla c(x,t)|^2\ dx-M\int_{\R^2}c(x,t)H(x)\ dx\ .\\
\end{eqnarray*}
Using inequality \fer{loginequality},
it is straightforward to prove that
\begin{equation*}
\begin{split}
-\f M{2\pi(1-\delta)}\iint_{\R^2\times \R^2}H(x)\log|x-y|n(y,t)\ dx dy
\ge -\f M{4\pi(1-\delta)}\int_{\R^2}n(x,t)\log(1+|x|^2)\ dx \\ -\f1{1-\delta} C(M)\ .
\end{split}
\end{equation*}
Finally, applying the logarithmic HLS inequality \fer{eq:logHLS} to $n$, we arrive exactly to the estimate
\begin{eqnarray*}
\mathcal E(t)&\ge&\left(1-\f M{8\pi(1-\delta)}\right)\int_{\R^2}n(x,t)\log n(x,t)\ dx
+\f M{8\pi(1-\delta)}\int_{\R^2}n(x,t)\log H(x)\ dx \\
&\ &\quad\quad\quad\quad\quad\quad\quad\quad\quad\quad\quad\quad
+\f12 \delta \int_{\R^2}|\nabla c(x,t)|^2\ dx-M\int_{\R^2}c(x,t)H(x)\ dx-\f 1{1-\delta} C(M)
\end{eqnarray*}
 and the proof follows as in section   \ref{sec:Onofri}.

\end{proof}

\begin{rem} Again, we see here that the mass $M$ cannot equal the critical value $8\pi$ (see Remark~\ref{rem:8pi}).
\end{rem}

\section{Regularizing effect}
\label{Sec:regularizing}
In the previous section it has been proved implicitly that under the hypothesis of Theorem \ref{th:aprioriestimate}, in particular for sub-critical mass $M$, the solution $n$ of system \fer{KS} is locally in time equi-integrable, i.e. there exists a modulus of equi-integrability $\omega(T;k)$, $T>0$ and $k>0$, such that
\begin{equation}
\sup_{0\le t\le T}\quad \int_{\R^2} \big(n(x,t)-k\big)_+\ dx \leq \omega(T;k) \quad \mbox{and}\quad \lim_{k\to+\infty}\omega(T;k) = 0\ .
\label{eq:equi-integrability}
\end{equation}
Indeed, obviously $\int_{\R^2} \big(n(x,t)-k\big)_+\ dx\le M$ for any $k>0$, while for $k>1$ we have
\begin{eqnarray}
\int_{\R^2} \big(n(x,t)-k\big)_+dx &\leq&
\f1{\log k}\int_{\R^2}\big(n(x,t)-k\big)_+\log n(x,t)\ dx\nonumber\\
&\leq&\f1{\log k}\int_{\R^2}n(x,t)\big(\log n(x,t)\big)_+dx
\le \f{C(1+t)}{\log k}\ .
\nonumber
\end{eqnarray}
\

In this section, following a now classical idea initiated in \cite{JL}, we will obtain {\em a priori} estimates for the $L^p$-norm of $n$, with the help of the equi-integrability property \fer{eq:equi-integrability}  and of the fact proved in Theorem \ref{th:aprioriestimate} that $\partial_tc\in L^2(0,T;L^2(\R^2))$. Since the hypothesis $\|n_0\|_{L^p(\R^2)}<\infty$ is not required, the following result is an hypercontractivity type result.

\begin{thm}\label{th:Lpestimates}
Let $T>0$ and $1<p<\infty$. Under the hypothesis of Theorem  \ref{th:aprioriestimate}, there exists a constant $C(T)$ not depending on $\| n_0\|_{L^p(\R^2)}$ such that
\begin{equation}
\int_{\R^2} n^p(x,t)\ dx \leq C(T) \left(1+t^{1-p}\right)\ ,\quad\forall\ 0<t\leq T \ ,
\label{eq:Lpestimates}
\end{equation}
i.e. the cell density $n(\cdot,t)$ belongs to $L^p(\R^2)$ for any positive time $t$.
\end{thm}
\begin{proof} Let $k>0$ to be chosen later.  We derive a non-linear differential inequality for the quantity $\displaystyle Y_p(t) := \int_{\R^2}(n(x,t)-k)_+^p\ dx$, which guarantees that the $L^p$-norm of $n$ remains finite whatever $\|n_0\|_{L^p(\R^2)}$ is (possibly infinite).
\\ 

{\sl First step : the differential inequality.} Multiplying the equation on $n$ in \fer{KS} by $p(n-k)_+^{p-1}$ yields, after integration by parts,
\begin{equation}
\begin{split}
\frac{d}{dt}\int _{\R^2} (n-k)_+^{p} \,dx = 
-4\f{(p-1)}p\int_{\R^2}  \big\vert\nabla(n-k)_+^{p/2}\big\vert^2 \,dx
-(p-1)\int_{\R^2}  (n-k)_+^{p}\Delta c \,dx \\ -pk\int_{\R^2}  (n-k)_+^{p-1}\Delta c \,dx\ .
\label{eq:Lpkcontrol}
\end{split} 
\end{equation}
There is some subtlety hidden here because we cannot use directly $-\Delta c=n-\alpha c$, as it is the case for the parabolic-elliptic system ($\e=0$). However, using the equation on $c$, one obtains
\begin{eqnarray}
\frac{d}{dt}\int _{\R^2}  (n-k)_+^{p} \,dx &\leq& 
-4\f{(p-1)}p\int_{\R^2}   \vert\nabla(n-k)_+^{p/2}\vert^2 \,dx \nonumber\\
&&+(p-1)\int_{\R^2}   (n-k)_+^{p+1}\,dx+(2p-1)k\int_{\R^2}   (n-k)_+^{p}\,dx
+pk^2\int_{\R^2}   (n-k)_+^{p-1}\,dx \nonumber\\
&&-\e(p-1)\int_{\R^2}   (n-k)_+^{p}\partial_tc\ dx
-\e pk\int _{\R^2}  (n-k)_+^{p-1}\partial_tc\ dx \ 
\label{eq:Lpkcontrol2}
\end{eqnarray}
and the additional non-linear terms $\int_{\R^2}   (n-k)_+^{p}\partial_tc\ dx$ and $\int _{\R^2}  (n-k)_+^{p-1}\partial_tc\ dx$ can be estimated in the following way. 
Using the Gagliardo-Nirenberg-Sobolev inequality
\[ 
\int _{\R^2} u^4(x)\ dx \leq C\int _{\R^2} u^2(x)\ dx \int _{\R^2} |\nabla u(x)|^2\ dx\ ,
\]
with $u = (n-k)_+^{p/2}$, we obtain
\begin{eqnarray} 
\Big|\int_{\R^2}  (n-k)_+^p \partial_t c\ dx  \Big|&\leq &  \Big(\int_{\R^2}  (n-k)_+^{2p}\ dx \Big)^{1/2} \|\partial_t c\|_{L^2(\R^2)}\nonumber\\
&\leq& C\Big(\int_{\R^2} (n-k)_+^{p}\ dx \Big)^{1/2} \Big(\int_{\R^2} \big|\nabla (n-k)_+^{p/2}\big|^2\ dx \Big)^{1/2} \|\partial_t c\|_{L^2(\R^2)}\nonumber\\
&\leq &  \e C(p)\|\partial_t c\|^2_{L^2(\R^2)}\int_{\R^2} (n-k)_+^{p} \ dx +  \frac2{\e p} \int_{\R^2}  \big|\nabla (n-k)_+^{p/2}\big|^2\ dx .
\label{estimatep}
\end{eqnarray}
On the other hand, by interpolation and the same Gagliardo-Nirenberg-Sobolev inequality as above, we have for $p\ge\f32$
\begin{eqnarray} 
\Big|\int_{\R^2}  (n-k)_+^{p-1} \partial_t c\ dx  \Big|&\leq &
\Big(\int_{\R^2}  (n-k)_+^{2(p-1)}\ dx \Big)^{1/2} \|\partial_t c\|_{L^2(\R^2)}\nonumber\\
&\leq &\Big(C(M,p)+C(p)\int_{\R^2}  (n-k)_+^{2p}\ dx \Big)^{1/2} \|\partial_t c\|_{L^2(\R^2)}\nonumber\\
&\leq &C(M,p)\|\partial_t c\|_{L^2(\R^2)}
+\e C(p,k)\|\partial_t c\|^2_{L^2(\R^2)}\int_{\R^2} (n-k)_+^{p} \ dx\nonumber\\ 
&\ &\quad\quad\quad\quad\quad\quad\quad\quad
+  \frac{(p-1)}{\e p^2 k} \int_{\R^2}  \big|\nabla (n-k)_+^{p/2}\big|^2\ dx .
\label{estimatep-1}
\end{eqnarray}
Inserting \fer{estimatep} and \fer{estimatep-1} in \fer{eq:Lpkcontrol2} gives for $p\ge\f32$
\begin{eqnarray}
\frac{d}{dt}\int _{\R^2}  (n-k)_+^{p} \,dx &\leq& 
-\f{(p-1)}p\int_{\R^2}   \vert\nabla(n-k)_+^{p/2}\vert^2 \,dx \nonumber\\
&&+(p-1)\int_{\R^2}   (n-k)_+^{p+1}\,dx+(2p-1)k\int_{\R^2}   (n-k)_+^{p}\,dx
+pk^2\int_{\R^2}   (n-k)_+^{p-1}\,dx \nonumber\\
&&+\e C\|\partial_t c\|^2_{L^2(\R^2)}\int_{\R^2} (n-k)_+^{p} \ dx 
+C\|\partial_t c\|_{L^2(\R^2)}\ .
\label{eq:Lpkcontrol3}
\end{eqnarray}
Next, we estimate the non-linear and negative contribution $-\f{(p-1)}p\int_{\R^2}   \vert\nabla(n-k)_+^{p/2}\vert^2 \,dx$ in term of $\int_{\R^2}   (n-k)_+^{p+1}\,dx$ and of the modulus of equi-integrability $\omega(T;k)$, with the help of the Sobolev's inequality. Indeed,
\begin{eqnarray}
\int_{\R^2}   (n-k)_+^{p+1}\,dx&=&\int_{\R^2}   \left((n-k)_+^{\f{(p+1)}2}\right)^2\,dx
\le C\left(\int_{\R^2}   \left\vert\nabla(n-k)_+^{\f{(p+1)}2}\right\vert \,dx\right)^2 \nonumber\\
&=&C(p)\left(\int_{\R^2} (n-k)_+^{\f12}  \vert\nabla(n-k)_+^{\f p2}\vert \,dx\right)^2
\nonumber\\
&\le&C(p)\int_{\R^2}   (n-k)_+\ dx\ \int_{\R^2}   \vert\nabla(n-k)_+^{p/2}\vert^2 \,dx
\nonumber\\
&\le&C(p)\omega(T;k)\int_{\R^2}   \vert\nabla(n-k)_+^{p/2}\vert^2 \,dx\ ,\quad\forall\ 0<t\le T\ .
\label{eq:estimategrad}
\end{eqnarray}
Moreover, since for $p\ge2$ it holds true that
\begin{equation}
\int_{\R^2}   (n-k)_+^{p-1}\,dx
\le \int_{\R^2}   (n-k)_+\,dx+\int_{\R^2}   (n-k)_+^{p}\,dx\ ,
\label{eq:estimatep-1}
\end{equation}
inserting \fer{eq:estimategrad} and \fer{eq:estimatep-1} in \fer{eq:Lpkcontrol3} gives for $p\ge2$ and $0<t\le T$
\begin{eqnarray}
\frac{d}{dt}\int _{\R^2}  (n-k)_+^{p} \,dx &\leq& (p-1)\left(1-\f1{pC(p)\omega(T;k)}\right)\int_{\R^2}   (n-k)_+^{p+1}\,dx
\nonumber\\
&\ &+C(1+\e\|\partial_t c\|^2_{L^2(\R^2)})\int_{\R^2}   (n-k)_+^{p}\,dx
+C\|\partial_t c\|_{L^2(\R^2)}+pk^2M\ .
\label{eq:Lpkcontrol4}
\end{eqnarray}
Finally, for any fixed $p$ we choose $k=k(p,T)$ sufficiently large such that
\begin{equation}
\delta:=\f1{pC(p)\omega(T;k(p,T))}-1>0\ .
\label{delta}
\end{equation}
This is clearly possible because $\omega(T;k)\to0$ as $k\to+\infty$. For this $k$ and using the interpolation
\begin{eqnarray*}
\int _{\R^2}(n-k)_+^{p} \,dx &\le&\left(\int _{\R^2}(n-k)_+\,dx\right)^{\f1p}
\left(\int _{\R^2}(n-k)_+^{p+1} \,dx\right)^{(1-\f1p)} \\ 
&\le& M^{\f1p}\left(\int _{\R^2}(n-k)_+^{p+1} \,dx\right)^{(1-\f1p)}\ ,
\end{eqnarray*}
we end up with the following differential inequality for $Y_p(t)$, $p\ge2$ fixed and $0<t\le T$
\begin{equation}
\begin{split}
 \dfrac{d}{dt} Y_p (t)\leq -(p-1)M^{-\f1{p-1}}\,\delta\  Y_p^\beta(t) + C_1\left(1 +\e \|\partial_t c(t)\|_{L^2(\R^2)}^2\right) Y_p(t) \\ + C_2\left(1 + \e\|\partial_t c(t)\|_{L^2(\R^2)}^2\right)\ ,
 \label{eq:Yp}
\end{split}
\end{equation}
where $\beta=\f p{p-1}>1$\ .
\\

{\sl Second step : estimate on $Y_p$, $p\ge2$.} Let us write the differential inequality \fer{eq:Yp} as follows for simplicity
\begin{equation}
 \dfrac{d}{dt} Y_p (t)\leq -\gamma\  Y_p^\beta(t) + f(t) Y_p(t) + f(t)\ ,\quad 0<t\le T\ ,
 \label{eq:Yp2}
\end{equation}
where $\gamma= (p-1)M^{-\f1{p-1}}\,\delta>0$ and $f(t)=\overline C\left(1 + \e\|\partial_t c(t)\|_{L^2(\R^2)}^2\right)$ with $\overline C=\max\{C_1,C_2\}$. Next we show that there exists a constant $C(T)$  not depending on $Y_p(0)$ such that
\begin{equation}
Y_p(t)\le C(T)\f1{t^{p-1}}\ ,\quad 0<t\le T\ ,
\label{eq:Ypstimate}
\end{equation}
by comparison of $Y_p(t)$ with positive solutions of the differential equation
\begin{equation}
\dfrac{d}{dt} Z_p (t)= -\gamma\  Z_p^\beta(t) + f(t) Z_p(t) + f(t)\ ,\quad 0<t\le T\ .
\label{eq:Zp}
\end{equation}
To do that, let $Z_p$ be a positive solution of \fer{eq:Zp} with $Z_p(0)\ge\left(\f{\overline C}{\gamma}\right)^{1/(\beta-1)}$ (such a solution exists from Carath\'eodory regularity of ordinary differential equations). Because $f(t)$ is positive,  $Z_p$ satisfies the differential inequality
\begin{equation}
\dfrac{d}{dt} Z_p (t)\ge -\gamma\  Z_p^\beta(t) + \overline C Z_p(t)\ ,\quad 0<t\le T\ .
\label{eq:Zpinequality}
\end{equation}
Therefore, since $\left(\f{\overline C}{\gamma}\right)^{1/(\beta-1)}$ is a constant solution of $Z'(t)=-\gamma\  Z^\beta(t) + \overline C Z(t)$, by comparison in \eqref{eq:Zpinequality} we get $Z_p(t)\ge \left(\f{\overline C}{\gamma}\right)^{1/(\beta-1)}$, $\forall\ 0\le t\le T$. As a consequence, from \fer{eq:Zp} $Z_p$ satisfies also the differential inequality
\begin{equation}
\dfrac{d}{dt} Z_p (t)\le -\gamma\  Z_p^\beta(t) + h(t) Z_p(t)\ ,\quad 0\le t\le T\ ,
\label{eq:Zpinequality2}
\end{equation}
where $h(t)=\left(1+\left(\gamma{\overline C}^{-1}\right)^{1/(\beta-1)}\right)f(t)$. Integrating \fer{eq:Zpinequality2} over $(0,t)$, it is straightforward to prove that
\begin{equation}
Z_p^{1-\beta}(t)\ge Z_p^{1-\beta}(0)e^{(1-\beta)\int_0^th(\tau)d\tau}+
\gamma(\beta-1)\int_0^te^{(1-\beta)\int_s^th(\tau)d\tau}ds\ ,\quad 0\le t\le T\ .
\nonumber
\end{equation}
Then, thanks to the estimate (vi) in Theorem \ref{th:aprioriestimate} and by the definition of $h(t)$, it holds true that
\begin{eqnarray*}
Z_p^{1-\beta}(t)&\ge& Z_p^{1-\beta}(0)e^{C(1-\beta)\e\int_0^T\|\partial_t c(s)\|_{L^2(\R^2)}^2ds}\ 
e^{C(1-\beta)t} \\ && \qquad +\f\gamma{C}\ e^{C(1-\beta)\e\int_0^T\|\partial_t c(s)\|_{L^2(\R^2)}^2ds}
(1-e^{C(1-\beta)t})\\
&\ge&\gamma(\beta-1)\ e^{C(1-\beta)(1+T)}e^{C(1-\beta)T}\ t=C(T)\ t\ ,\quad 0\le t\le T\ ,
\end{eqnarray*}
and therefore
\begin{equation}
Z_p(t)\le C(T)\f1{t^{p-1}}\ ,\quad 0\le t\le T\ .
\label{eq:Zpinequality3}
\end{equation}
where $C(T)$ does not depend on $Z_p(0)$. 

Finally, let $Z_p$ be a positive solution of \fer{eq:Zp} with $Z_p(0)=\max\left\{\left(\f{\overline C}{\gamma}\right)^{1/(\beta-1)},Y_p(0)\right\}$. Again by comparison, $Y_p(t)\le Z_p(t)$, $\forall\ 0<t\le T$ and \fer{eq:Ypstimate} follows by \fer{eq:Zpinequality3}.
\\

{\sl Third step : $L^p$ regularity of $n$.} To conclude, it is sufficient to observe that for any $k>0$ we have
\begin{eqnarray}
\int_{\R^2}n^p(x,t)\ dx&=&\int_{\{n\le 2k  \}}n^p(x,t)\ dx+\int_{\{n>2k  \}}n^p(x,t)\ dx\nonumber\\ 
&\le&(2k  )^{p-1}M+ 2 ^p\int_{\{n>2k  \}}(n(x,t)-k)^p\ dx\nonumber \\
&\le&(2k  )^{p-1}M+2^p\int_{\R^2}(n(x,t)-k)_+^p\ dx\ ,
\label{np}
\end{eqnarray}
where the inequality $x^p\le2^p(x-k)^p$, for $x\ge k 2 $, has been used. Therefore, estimate \fer{eq:Lpestimates} follows for any $p\ge2$ by \fer{eq:Ypstimate} and \fer{np} choosing $k=k(p,T)$ sufficiently large such that \fer{delta} holds true. For $1<p<2$, the Theorem follows by interpolation.

\end{proof}

\section{Global Existence}
\label{Sec:Globalexistence}
We are now able to prove Theorem \ref{the:existence} collecting all the estimates proved in the previous sections. In order to do that, we need first to regularize the chemotaxis  system \fer{KS} and then to prove that the {\em a priori} estimates hold true and pass to the limit. Being this procedure quite technical and usual, we just sketch the proof. For the parabolic-elliptic case with $\alpha=0$ one can consult for example \cite{BDP}, where the regularizing procedure has been written in full details.

The regularized system that we consider is
\begin{equation}
\left\{\begin{array}{rcl}
\dfrac{\partial n^\sigma}{\partial t}  & =& \Delta n^\sigma - \nabla\cdot( n^\sigma\nabla c^\sigma)\ ,\quad t>0,\ x\in \R^2,\vspace{.2cm}\\
\e\dfrac{\partial c^\sigma}{\partial t} &= &\Delta c^\sigma + n^\sigma*\rho^\sigma - \alpha c^\sigma\ , \quad t>0,\ x\in \R^2,\vspace{.2cm}\\
n^\sigma(\cdot,0)&=&n_0*\rho^\sigma ,\quad c^\sigma(\cdot,0)=c_0*\rho^\sigma\ ,\quad\quad x\in \R^2\ ,
\end{array}
\right. 
\label{KSregularized}
\end{equation}
for some regularizing kernel $\rho^\sigma(x)= \f{1}{\sigma^2}\rho(\f{x}{\sigma})$ with $\rho\in \cad^+(\R^2)$ and $\int_{\R^2} \rho(x)\ dx=1$. The first step is to prove the global existence of a smooth solution $(n^\sigma,c^\sigma)$ of \fer{KSregularized}. This result  can be obtained through a Picard fixed-point method, writing $(n^\sigma,c^\sigma)$ as
$$
n^\sigma(t)=G(t)* (n_0*\rho^\sigma)-\int_0^t \nabla G(t-s)*(n^\sigma(s)\nabla c^\sigma(s))\ ds\ ,
$$
$$
c^\sigma(t)=e^{-\alpha t/\e} G(t/\e)* (c_0*\rho^\sigma) +\frac1\e\int_0^{t/\e} e^{\alpha(s-t)/\e} G((t-s)/\e)* (n^\sigma(s/\e)*\rho^\sigma) \ ds\ ,
$$
with $G(x,t)=\frac{1}{4\pi t}e^{-|x|^2/(4t)}$ the heat kernel in $\R^2$, because the nonlinearity is Lipschitz (see \cite{B98}). The smoothness of $(n^\sigma,c^\sigma)$ follows by the regularizing property of the heat equation and by the smoothness of the initial data. The non-negativity of $(n^\sigma,c^\sigma)$ follows by the maximum principle. As a consequence, the {\em a priori} estimates in Theorem \ref{th:aprioriestimate} hold true for the regularized solution $(n^\sigma,c^\sigma)$. 
These {\em a priori} bounds give a global in time control of the energy and energy dissipation, still uniform in $\sigma$. This gives space compactness on $c$ and $\nabla n$. The Lions-Aubin compactness method gives the required time compactness.

\section{Blow-up}
\label{Sec:blowup}
In this section we shall prove Theorem \ref{the:blowup} which states a blow-up result for the parabolic-elliptic system \fer{KS} and for super-critical mass, i.e. $M>8\pi$. Since the same result in the case $\e=\alpha=0$ has been given in \cite{BDP}, we assume $\alpha>0$ so that  $c=B_\alpha*n$. The result in \cite{BDP} will be obtained in the limit $\alpha\to0$.

Let us denote $I(t)=\int_{\R^2}|x|^2n(x,t)\ dx$. Then, $I(t)$ satisfies the differential equation
$$
\f d{dt}I(t)=4M+2\int_{\R^2}n(x,t)\ x\cdot(\nabla B_\alpha*n(t))(x)\ dt\ .
$$
Using the computation
$$
\nabla B_\alpha (z)=-\f z{8\pi}\int_0^{+\infty}\f1{t^2}\ e^{-\f{|z|^2}{4t}-\alpha t}\ dt
=-\f1{2\pi}\f z{|z|^2}\int_0^{+\infty}e^{-s-\alpha\f{|z|^2}{4s}}\ ds
$$
and denoting $g_\alpha(z)=\int_0^{+\infty}e^{-s-\alpha\f{|z|^2}{4s}}\ ds$, we obtain 
\begin{eqnarray}
\f d{dt}I(t)&=&4M-\f1\pi\iint_{\R^2\times\R^2}n(x,t)\f{x\cdot(x-y)}{|x-y|^2}g_\alpha(x-y)n(y,t)\ dydx\nonumber\\ 
&=&4M-\f1{2\pi}\iint_{\R^2\times\R^2}n(x,t)g_\alpha(x-y)n(y,t)\ dydx\ .
\label{momenteq}
\end{eqnarray}
Since $g_\alpha$ is a positive radial decreasing function such that $g_\alpha(z)\le1$ and $g_\alpha(z)\to1$ as $\alpha\to0$ for all $z\in\R^2$, we recover from \fer{momenteq} the result in \cite{BDP}. For $\alpha>0$ we have
\begin{equation}
\f d{dt}I(t)=4M\left(1-\f M{8\pi}\right)+\f1{2\pi}\int_{\R^2}\int_{\R^2}n(x,t)[1-g_\alpha(x-y)]n(y,t)\ dydx
\label{2momenteq}
\end{equation}
and one has to estimate the second term in the right hand side of \fer{2momenteq} in term of $I(t)$. Denoting $r=|z|$, we observe that
\begin{equation}
\f d{dr}(1-g_\alpha(r)) =\f\alpha2\ r\int_0^{+\infty}\f1se^{-s-\alpha\f{r^2}{4s}}\ ds
=2\pi\alpha rB_1(\sqrt\alpha r)\ .
\label{galpha}
\end{equation}
Then, reasoning as in Lemma \ref{lem:Besselinequality}, we easily deduce from \fer{galpha} that
\begin{equation}
\f d{dr}(1-g_\alpha(r))\le\sqrt\alpha\ K\ ,\quad\forall\ 0<r<\f1{\sqrt\alpha}\ ,
\label{galpha2}
\end{equation} 
where $K= 2\pi\sup_{\rho\in (0,1)}\rho B_1(\rho)<+\infty$. As a consequence of \fer{galpha2} and since $0\le1-g_\alpha(z)\le1$ for all $z\in\R^2$, we obtain the estimate
\begin{equation}
0\le1-g_\alpha(z)\le \sqrt{\alpha}\ \mathcal C|z|\ ,\quad\forall z\in\R^2\ ,
\label{galpha3}
\end{equation} 
with $\mathcal C=\max\{K, 1 \}$.

Finally, we insert \fer{galpha3} into \eqref{2momenteq} and we have
\begin{eqnarray*}
 \f d{dt}I(t) &\leq& 4M\left(1-\f M{8\pi}\right)+\f{\sqrt\alpha}{2\pi} \mathcal C\iint_{\R^2\times\R^2}n(x,t)|x-y|n(y,t)\ dydx \\
&\leq & 4M \left(1-\f M{8\pi}\right) + \f{\sqrt\alpha}{\pi} \mathcal CM \int_{\R^2} |x|n(x,t)\ dx \\
&\leq & 4M \left(1-\f M{8\pi}\right) + \f{\sqrt\alpha}{\pi} \mathcal CM^{3/2} \sqrt{I(t)} \ ,
\end{eqnarray*}
i.e.
\begin{equation}
I(t)\le I(0)+\int_0^tf(I(s))\ ds\ ,
\label{eq:I(t)}
\end{equation} 
where $f(\lambda)=\f M{2\pi}(8\pi-M)+\f{\sqrt\alpha}{\pi} \mathcal CM^{3/2}\lambda^{\f12}$. Since $f$ is a strictly increasing function such that $f(\lambda^*)=0$ for $\lambda^*=\f1{4\alpha\mathcal C^2M}(M-8\pi)^2$, the hypothesis \fer{I(0)} gives us $I(0)<\lambda^*$ and $f(I(0))<0$. Therefore, $\int_0^tf(I(s))\ ds<0$ as soon as \fer{eq:I(t)} holds true and $I(t)\le I(0)+tf(I(0))$ as soon as \fer{eq:I(t)} holds true, i.e. the second momentum becomes necessarily non-positive for $t\ge-\f{I(0)}{f(I(0))}=\f{2\pi I_0}{M(M-8\pi-2\mathcal C\sqrt{\alpha MI_0})}$ expressing in this way the formation of a singularity before.
\section{Appendix} 
\label{sec:appendix}
\subsection{The entropy minimization Lemma} 
\label{ap:entropy}
The entropy minimization Lemma \ref{lem:entropy} is very classical and several proofs are available. For example, we may use the Jensen's inequality with respect to the probability measure $\dfrac nM dx$, 
and we deduce
%
\begin{eqnarray}
\f1M\int_{\R^2}(n(x)\log n(x)-n(x)\psi(x))\ dx&=&-\int_{\R^2}\log\left(\f{e^{\psi(x)}}{n(x)}\right)\ \frac{n(x)}{M}dx
\nonumber\\
&\ge&-\log\left(\f1M\int_{\R^2}e^{\psi(x)}\ dx\right)=\f1M E(\overline n;\psi)\ .
\nonumber
\end{eqnarray}
One can also appeal to the Legendre transform of the functional $\int_{\R^2} n(x)\log n(x)\ dx$ to find out \fer{entropyequality}. However, in these ways we loose the identity in  \fer{entropyequality}. Therefore, the following proof, in the line of  \cite{CJMTU}, is  more complete.
\\

\begin{proof}[Proof of the entropy minimization Lemma \ref{lem:entropy}] First of all, by the definition of $\overline{n}$ we have that $\overline n\in{\cal U}$ and 
\begin{equation}
\log \overline{n} = \psi + \log\left(M\Big/\int_{\R^2} e^{\psi}\ dx\right).
\label{lognbar}
\end{equation} 
Therefore, the entropy functional $E$ is finite in $\overline n$ and it takes the value $E(\overline n;\psi)=M\log\left(M/\int_{\R^2} e^{\psi}\ dx\right)$. Next, it is easy to see that for any $n$ in ${\cal U}$, $E(n;\psi)$ and $RE(n|\overline{n})$ are finite or infinite in the same time and that \fer{entropyequality} holds true. Indeed, from \fer{lognbar} we deduce $n\log( n/{\overline n})=n\log n-n\psi-n\log\left(M/\int_{\R^2} e^{\psi}\ dx\right)$. The non-negativity of $RE(n|\overline{n})$ over the set ${\cal U}$ follows by the computation
$$
RE(n|\overline{n})=\int_{\R^2}\ (n(x)\log n(x)-\overline n(x)\log\overline n(x)-(\log\overline n(x)+1)(n(x)-\overline n(x)))\ dx \ ,
$$
and the convexity of the function $u\mapsto u\log u$.
\end{proof} 
\subsection{The chemical energy minimization Lemma} 
\label{ap:chemicalminimization}
\begin{proof}[Proof of the chemical energy minimization Lemma \ref{lem:chemicalmin}]
Let us start first with the basic regularity properties of $\overline c$ (see\cite{LiebLoss}). If $\alpha>0$, then $\overline c\in L^p(\R^2)$, for all $p\in[1,\infty)$ since $f\in L^1(\R^2)$ by hypothesis. On the other hand, when $\alpha=0$, $\overline c\in L_{loc}^1(\R^2)$ since $f\in L^1(\R^2)\cap L^1(\R^2,\log(1+|x|^2)dx)$. Moreover \fer{gradc} holds true in ${\cal D}'(\R^2)$ as well as $-\Delta\overline c+\alpha\overline c=f$ for $\alpha\ge0$.

Next we aim to justify the integration by parts arising in the computation \fer{bo}. Let us assume that  we know {\em a priorily} that $|\nabla\overline c|\in L^2(\R^2)$. Then $c\Delta\overline c\in L^1(\R^2)$ for all $c\in H^1(\R^2)$ and the following partial integration holds true
\begin{equation}
\int_{\R^2}\nabla c(x)\cdot\nabla\overline c(x)\ dx=-\int_{\R^2}c(x)\Delta \overline c(x)\ dx=\int_{\R^2}c(x)(f(x)-\alpha\overline c(x))\ dx\ .
\label{partialintegr}
\end{equation} 
Indeed, thanks to the hypotesis and the basic properties of $\overline c$ we can decompose $-\Delta\overline c=g_1+g_2$ with $g_1\ge0$ in $L^1_{loc}(\R^2)$ defined as 
$$
g_1:=f+\alpha(\overline c)_-\quad\hbox{if }\alpha>0
\quad\hbox{and}\quad
g_1:=f_1\quad\hbox{if }\alpha=0
$$  
and $g_2\in L^2(\R^2)$ defined as
$$
g_2:=-\alpha(\overline c)_+\quad\hbox{if }\alpha>0
\quad\hbox{and}\quad
g_2:=f_2\quad\hbox{if }\alpha=0\ ,
$$
(see \cite[Theorem 7.7]{LiebLoss}). Moreover, $F_\alpha(c;f)$ is finite. 

As a consequence of \fer{partialintegr}, we obtain easily \fer{bo} since
\begin{eqnarray*}
\f12\int_{\R^2}|\nabla c-\nabla\overline c|^2+\f\alpha2\int_{\R^2}(c-\overline c)^2\ dx &= &
F_\alpha(c;f)+\f12\int_{\R^2}|\nabla\overline c|^2\ dx+\f\alpha2\int_{\R^2}\overline c^2\ dx 
\\
&=&F_\alpha(c;f)-F_\alpha(\overline c;f)\ .
\end{eqnarray*}
\

It remains to prove that $|\nabla \overline c|\in L^2(\R^2)$. Let us consider first the case $\alpha=0$ that is much more complicated then the case $\alpha>0$ since the fundamental solution $E_2$ does not lie in any $L^p(\R^2)$ spaces and $\nabla\overline c=\nabla E_2*f$ in general does not lie in $L^2(\R^2)$ because of the critical fractional Sobolev embedding (see \cite{LiebLoss}). The case $\alpha>0$ is considered subsequently.

In the case $\alpha=0$, the hypothesis $\int_{\R^2}f(x)\ dx=0$ allows us to substract any function of $x$ to the integrand inside $\nabla\overline c$. Let us consider two radii $0<r<1<R$ with $R>e-1$ and let us denote ${\cal C}B(0,r)=\R^2\setminus B(0,r)$. Then, we split $\nabla\overline c$ in the following way
\begin{eqnarray}
\nabla\overline c(x)&=&-\f1{2\pi}\int_{\R^2}\left(\f{x-y}{|x-y|^2}-\f{x}{1+|x|^2}\ \bone _{{\cal C}B(0,R)}(x)\right)f(y)\ dy
\nonumber\\
&=&-\f1{2\pi}\int_{\R^2}\left(\f{x-y}{|x-y|^2}-\f{x}{1+|x|^2}\ \bone _{{\cal C}B(0,R)}(x)\right)\ \bone _{B(0,r)}(x-y)f(y)\ dy
\nonumber\\
& &-\f1{2\pi}\int_{\R^2}\left(\f{x-y}{|x-y|^2}-\f{x}{1+|x|^2}\ \bone _{{\cal C}B(0,R)}(x)\right)\ \bone _{{\cal C}B(0,r)}(x-y)f(y)\ dy\nonumber\\
&=:& I_1(x)+I_2(x)
\end{eqnarray}
and we will show separately that $|I_1|$ and $|I_2|$ lie in $L^2(\R^2)$. 

Concerning $I_1$ we have for all $x\in\R^2$,
\begin{align}
|I_1(x)|\le & \f1{2\pi}\int_{\R^2}\f1{|x-y|}\ \bone _{B(0,r)}(x-y)|f(y)|\ dy
\nonumber \\ &\qquad +\f{|x|\bone _{{\cal C}B(0,R)}(x)}{2\pi(1+|x|^2)}\ \int_{\R^2}\bone _{B(0,r)}(x-y)|f(y)|\ dy\ .
\label{eq:I1}
\end{align}
Denoting $\Omega_x=\left\{y\in\R^2\ :\ f_1(y)>\f1{|x-y|}\right\}$, the first integral in the right hand side of \fer{eq:I1} can be split again in the following way
\begin{eqnarray}
\int_{\R^2}\f1{|x-y|}\ \bone _{B(0,r)}(x-y)|f(y)|\ dy&\le&
\int_{\R^2}\f1{|x-y|}\ \bone _{B(0,r)}(x-y)f_1(y)\bone _{\Omega_x}(y)\ dy
\nonumber\\
& &+\int_{\R^2}\f1{|x-y|}\ \bone _{B(0,r)}(x-y)f_1(y)\bone _{{\cal C}\Omega_x}(y)\ dy
\nonumber\\
& &+\int_{\R^2}\f1{|x-y|}\ \bone _{B(0,r)}(x-y)|f_2(y)|\ dy
\nonumber\\
&=:&I_{1,1}(x)+I_{1,2}(x)+I_{1,3}(x)\ .
\label{eq:I1+I2}
\end{eqnarray}
Since for $y\in B(x,r)\cap\Omega_x$ it holds true that $f_1(y)>\f1{|x-y|}>\f1r>1$, we have
\begin{equation}
0\le I_{1,1}(x)\le\int_{\R^2}\f1{|x-y|}\f1{\log\left(\f1{|x-y|}\right)}\ \bone _{B(0,r)}(x-y)f_1(y)\log f_1(y)\bone _{\Omega_x}(y)\ dy\ .
\label{eq:I11}
\end{equation}
Therefore,  since the right hand side of \fer{eq:I11} belongs to $L^2(\R^2)$ thanks to the Young's inequality, $I_{1,1}\in L^2(\R^2)$ too. On the other hand, since for $y\in B(x,r)\cap{\cal C}\Omega_x$ it holds true that $f_1(y)\le\f1{|x-y|}$, we have
\begin{equation}
0\le I_{1,2}(x)\le\int_{\R^2}\f1{|x-y|^{3/2}}\ \bone _{B(0,r)}(x-y)\sqrt{f_1(y)}\bone _{{\cal C}\Omega_x}(y)\ dy\ .
\label{eq:I12}
\end{equation}
Again by the Young's inequality the right hand side of \fer{eq:I12} belongs to $L^2(\R^2)$ and so $I_{1,2}$ too. The term $I_{1,3}$ in \fer{eq:I1+I2} belongs to $L^2(\R^2)$ easily by Young's  inequality since $f_2\in L^2(\R^2)$. Finally, since for $x\in {\cal C}B(0,R)$ and $y\in B(x,r)$ we have $|y|\ge|x|-r\ge R-r>0$, the second term in the right hand side of \fer{eq:I1} can be dominated in the following way
\begin{align}
&\f{|x|\bone _{{\cal C}B(0,R)}(x)}{1+|x|^2}\ \int_{\R^2}\bone _{B(0,r)}(x-y)|f(y)|\ dy \nonumber \\
&\qquad \le\f{|x|\bone _{{\cal C}B(0,R)}(x)}{(1+|x|^2)\log(1+(|x|-r)^2)}\int_{\R^2}|f(y)|\log(1+|y|^2)\ dy\ .
\label{eq:I1second}
\end{align}
Therefore, it belongs to $L^2(\R^2)$. Collecting \fer{eq:I1}, \fer{eq:I1+I2}, \fer{eq:I11}, \fer{eq:I12} and \fer{eq:I1second}, we obtain that $|I_1|$ belongs to $L^2(\R^2)$. 

Concerning $I_2$, it is enough to prove that $|I_2|\in L^2({\cal C}B(0,R))$ since $|I_2|\in L^\infty(\R^2)$. Let $x\in{\cal C}B(0,R)$ and let us define $\Omega'_x=\left\{y\in\R^2\ :\ |y|<|x|/\log(1+|x|)\right\}$. Then,
\begin{eqnarray}
|I_2(x)|&\le& \f1{2\pi}\int_{\R^2}\left|\f{x-y}{|x-y|^2}-\f{x}{1+|x|^2}\right|\ \bone _{{\cal C}B(0,r)}(x-y)|f(y)|\bone _{\Omega'_x}(y)\ dy
\nonumber\\
& &\quad+\f1{2\pi}\int_{\R^2}\left|\f{x-y}{|x-y|^2}-\f{x}{1+|x|^2}\right|\ \bone _{{\cal C}B(0,r)}(x-y)|f(y)|\bone _{{\cal C}\Omega'_x}(y)\ dy
\nonumber\\
&=:&I_{2,1}(x)+I_{2,2}(x)\ .
\label{est1}
\end{eqnarray}
For $x\in{\cal C}B(0,R)$ and $y\in\Omega'_x$ it holds true that $|x|-|y|>|x|\left(1-\f1{\log(1+|x|)}\right)>0$ since $R>e-1$ and
\begin{eqnarray}
\left|\f{x-y}{|x-y|^2}-\f{x}{1+|x|^2}\right|&\le&
\f{|(x-y)(1+|x|^2)-x(|x|^2+|y|^2-2x\cdot y)|}{(|x|-|y|)^2(1+|x|^2)}
\nonumber\\
&\le&\f{|x|+|y|+3|y||x|^2+|x||y|^2}{(1+|x|^2)}\f{\log^2(1+|x|)}{|x|^2(\log(1+|x|)-1)^2} 
\nonumber
\label{eq:I21}
\\
&\le&C\f{\log(1+|x|)}{|x|(\log(1+|x|)-1)^2}=:h(x)\ .
\end{eqnarray}
Therefore, $|I_{2,1}(x)|\le \f1{2\pi}\norm{f}_{L^1(\R^2)}h(x)$, which implies that $I_{2,1}\in L^2({\cal C}B(0,R))$. 
Finally, 
\begin{eqnarray}
0\le I_{2,2}(x)&\le&\f1{2\pi}\int_{\R^2}\f1{|x-y|}\ \bone _{\mathcal C B(0,r)}(x-y)|f(y)|\bone _{{\cal C}\Omega'_x}(y)\ dy\nonumber\\
&& \qquad +\f{|x|}{2\pi(1+|x|^2)}\  \int_{\R^2}\bone _{\mathcal C B(0,r)}(x-y)|f(y)|\bone _{{\cal C}\Omega'_x}(y)\ dy\ .
\label{eq:I22}
\end{eqnarray}
The second integral in the right hand side of \fer{eq:I22} belongs to $L^2({\cal C}B(0,R))$ because $|y|\ge\f{|x|}{\log(1+|x|)}$ for $y\in{\cal C}\Omega'_x$ so that
\begin{eqnarray}
&&\f{|x|}{(1+|x|^2)}\int_{\R^2}\bone _{\mathcal C B(0,r)}(x-y)|f(y)|\bone _{{\cal C}\Omega'_x}(y)\ dy\qquad\nonumber\\
&&\qquad\le\dfrac {|x|}{(1+|x|^2) \log\left(1+\frac{|x|^2}{ \log^2(1+|x|)}\right)}\int_{\R^2}\log(1+|y|^2)|f(y)|\ dy\in L^2({\cal C}B(0,R)).
\label{est2}
\end{eqnarray}
In order to control the first integral in the right hand side of \fer{eq:I22}, we observe that for all $x\in {\cal C}B(0,R)$ and $y\in{\cal C}\Omega'_x$ it holds true that
\begin{eqnarray} 
\log\left(1+\frac{\frac12|x-y|}{ \log(1+\frac12|x-y|)}\right) &\! \leq\!\! &\max\left( \log\left(1+\frac{|x|}{ \log(1+|x|)}\right)\ ;\ \log\left(1+\frac{|y|}{\log(1+|y|)}\right)\right)
\label{firstest}\\ 
&  \leq& \log\left(1+|y|^2\right)\ .
\label{eq:lastbound}  
\end{eqnarray}
Indeed, the function $\varphi(s)=\f s{\log(1+s)}$ is strictly increasing for $s>-1$ such that $\varphi(e-1)=e-1$. Then using the inequality $\f12|x-y|\le\max(|x|;|y|)$, the first inequality \fer{firstest} follows. Moreover, for $x\in {\cal C}B(0,R)$ and $y\in{\cal C}\Omega'_x$ we have 
$$
e-1\le\f{|x|}{\log(1+|x|)}\le|y|\quad\hbox{and}\quad
\f{|y|}{\log(1+|y|)}\le|y|\le|y|^2
$$
and the second inequality \fer{eq:lastbound} follows too. Therefore,
\begin{eqnarray}
&&\int_{\R^2}\f1{|x-y|}\ \bone _{\mathcal C B(0,r)}(x-y)|f(y)|\bone _{{\cal C}\Omega'_x}(y)\ dy \nonumber\\
&\le & 
\int_{\R^2}  \dfrac1{|x-y|}   \dfrac1{ \log\left(1+\frac{\frac12|x-y|}{ \log(1+\frac12|x-y|)}\right)}\ \bone _{\mathcal C B(0,r)}(x-y)\log\left(1+|y|^2\right) |f(y)|\bone _{{\cal C}\Omega'_x}(y)\ dy
\label{eq:I22left} 
\end{eqnarray}
Collecting \fer{est1}, \fer{eq:I21}, \fer{eq:I22}, \fer{est2} and \fer{eq:I22left} we have proved that $|I_2|\in L^2({\cal C}B(0,R))$ and the case $\alpha=0$ is solved.

When $\alpha >0$, we have that $|\nabla B_\alpha(z)|=\f1{2\pi}\f1{|z|}g_\alpha(z)$ with $g_\alpha(z)=\int_0^{+\infty}e^{-s-\alpha\f{|z|^2}{4s}}\ ds$, (see Section \ref{Sec:blowup}). Therefore, $|\nabla B_\alpha(z)|$ has the same singularity as $|\nabla E_2(z)|$ in $z=0$ but $|\nabla B_\alpha(z)|\to 0$ exponentially as $|z|\to+\infty$. As a consequence, we can prove that $\nabla\overline c=\nabla B_\alpha*f\in L^2(\R^2)$ using the previous technique and without substracting any function of $x$ to the integrand in $\nabla\overline c$.
\end{proof} 
\subsection{A remark on the case $\e=\alpha=0$} 
\label{ap:eps=alpha=0} 
In the case of the parabolic-elliptic system \fer{KS} with $\e=\alpha=0$ and 
$$
c(x,t)=-\f1{2\pi}\int_{\R^2}\log|x-y|n(y,t)dy\ ,
$$
as far as global existence is concerned, it is sufficient to assume that the cell density $n_0$ satifies both 
\begin{equation} \int_{\R^2}n_0(x) \log n_0(x)\ dx<\infty \quad \mbox{and}\quad \int_{\R^2}n_0(x)\log(1+|x|^2)\ dx<\infty\ .
\label{eq:minimalHyp}\end{equation}
In fact, assumptions \eqref{eq:minimalHyp} are the optimal ones for several viewpoints. First, they are minimal for applying the logarithmic HLS inequality \fer{eq:logHLS}. Furthermore, the combination of these two ensures that the mass does not escape to infinity by Lemma \ref{lem:controlmass} and the balance is optimal again. Last but not least, it can be proved using free energy methods that conditions \eqref{eq:minimalHyp}  are indeed propagated along the solutions, with local in time bounds. \\
Following the lines of computation \eqref{eq:cH} for instance, we have to estimate 
$\int_{\R^2}c(x,t) H(x)\ dx$ from above. This can be done with the following calculation,
\begin{eqnarray*}
 \int_{\R^2}c(x,t) H(x)\ dx & = & -\frac1{2\pi} \iint_{\R^2\times\R^2} n(y,t) \log|x-y| H(x)\ dydx \\
& = & \frac 1{8\pi} \int_{\R^2}n(y,t) \log H(y) \ dy + C\ ,
\end{eqnarray*}
which can be viewed as an integration by parts knowing \eqref{eq:cancellation} (see also \cite{CarlenLoss} where the Euler-Lagrange equation for $H$ is clearly noticed). The latter is then trivially bounded from above because $\int_{\R^2}n(y,t) \log H(y) \ dy$ is nonpositive.

Concerning the blow-up of solutions, the assumption on the second momentum, $\int_{\R^2}|x|^2 n_0(x) \ dx<+\infty$, seems to be however crucial. For instance, for the critical mass value $M=8\pi$, there exists a family of stationary states with infinite second momentum but finite ``logarithmic" momentum, for which blow-up does not occur obviously (see \cite{Biler06},\cite{BCM}).

\subsection{The duality} 
\label{duality}
There is a true duality between the Onofri and the logarithmic Hardy-Littlewood-Sobolev inequalities (see \cite{CarlenLoss},\cite{Beckner}). We followed this idea all along this paper. In this appendix we give a formal proof of this duality in the whole space $\R^2$ for the sake of completeness, and we next derive another proof of Lemma \ref{lem:Besselinequality}.

First of all, we write the Onofri's  inequality as
\begin{equation}
\f1{8\pi}\log\left(\int_{\R^2}e^{8\pi u(x)}H(x)\ dx\right)\leq\int_{\R^2}u(x)H(x)\ dx+\f12\int_{\R^2}\vert\nabla u(x)\vert^2\ dx\ .
\label{eq:OnoR2bis}
\end{equation}
Let $f$ be a function satisfying the hypotheses of Lemma \ref{lem:logHLS} and $\int_{\R^2} f(x)\ dx=1$ without loss of generality. By the minimization procedure (see Lemma \ref{lem:chemicalmin} and Appendix \ref{ap:chemicalminimization}) we have  
\begin{equation*}
\begin{split}
&-\f12\int_{\R^2}(f-H)(x)(E_2*(f-H))(x)\ dx
\\
&\qquad  =\min_{u}\left\{\f12\int_{\R^2}|\nabla u(x)|^2\ dx -
\int_{\R^2}(f-H)(x)u(x)\ dx\right\}\\
& \qquad \geq  \min_{u}\left\{ \f1{8\pi}\log\left(\int_{\R^2}e^{8\pi u(x)}H(x)\ dx\right) -
\int_{\R^2}f(x)u(x)\ dx\right\} \\
&\qquad \geq \min_{u}\left\{ - \frac1{8\pi} \int_{\R^2} f(x)\log f(x) - f(x)\big( 8\pi u(x)+\log H(x) \big)\ dx-
\int_{\R^2}f(x)u(x)\ dx\right\} \\
 &\qquad \geq - \frac1{8\pi} \int_{\R^2} f(x)\log f(x)\ dx +  \frac 1{8\pi}\int_{\R^2} f(x) \log H(x)\ dx\ ,
\end{split}
\end{equation*}
from the entropy minimization Lemma \ref{lem:entropy} (see Appendix \ref{ap:entropy}). On the other hand, the Euler-Lagrange formula for $H$ writes \cite{CarlenLoss},
\[ - \int_{\R^2} \log|x-y|H(y)\ dy = \frac14\log H(x) + \frac14 \int_{\R^2} H(x)\log H(x)\ dx +C_0\ . \]
It is in fact the dual formulation of the cancellation property \eqref{eq:cancellation}. We deduce that 
\begin{equation*}
\begin{split}
&\f1{4\pi}\iint_{\R^2\times\R^2} f (x)\log |x-y|f(y)\ dxdy
+\f1{4\pi}\iint_{\R^2\times\R^2} H(x)\log |x-y|H(y)\ dxdy\\
&\qquad \geq 
 \frac1{2\pi} \iint_{\R^2\times\R^2} f(x)\log|x-y|H(y)\ dxdy - \frac1{8\pi} \int_{\R^2} f(x)\log f(x)\ dx \\
&
\qquad \qquad \qquad +  \frac 1{8\pi}\int_{\R^2} f(x) \log H(x)\ dx  \\
&\qquad\geq - \frac1{8\pi} \int_{\R^2} f(x)\log f(x)\ dx
 - \frac1{8\pi} \int_{\R^2} H(x)\log H(x)\ dx +\frac{C_0}{2\pi}\ . 
\end{split}
\end{equation*}
We have recovered the logarithmic Hardy-Littlewood-Sobolev inequality with the sharp constant.

Following the previous lines there is another way to derive the modified inequality for the kernel $B_\alpha$ (Lemma \ref{lem:Besselinequality}). 
We proceed as above, more directly because the Bessel kernel has nicer properties at infinity:
\begin{equation*}
\begin{split}
&-\f12\int_{\R^2}f(x)(B_\alpha*f)(x)\ dx
\\
&\qquad  =\min_{u}\left\{\f12\int_{\R^2}|\nabla u(x)|^2\ dx+\f\alpha2\int_{\R^2}u^2(x)\ dx -
\int_{\R^2}f(x)u(x)\ dx\right\}\\
& \qquad \geq  \min_{u}\left\{ \f1{8\pi}\log\left(\int_{\R^2}e^{8\pi u(x)}H(x)\ dx\right) -
\int_{\R^2}f(x)u(x)\ dx+\f\alpha2\int_{\R^2}u^2(x)\ dx - \int_{\R^2} u(x) H(x)\ dx\right\} \\
& \qquad \geq  \min_{u}\left\{ - \frac1{8\pi} \int_{\R^2} f(x)\log f(x) - f(x)\big( 8\pi u(x)+\log H(x) \big)\ dx  -
\int_{\R^2}f(x)u(x)\ dx \right. \\ 
&\qquad\qquad\qquad \left.+\f\alpha2\int_{\R^2}u^2(x)\ dx - \int_{\R^2} u(x) H(x)\ dx \right\} \\
& \qquad \geq - \frac1{8\pi} \int_{\R^2} f(x)\log f(x) \ dx + \frac1{8\pi}\int_{\R^2} f(x)\log H(x)\ dx + C(\alpha)\ ,
\end{split}
\end{equation*}
because we have for $\alpha>0$:
\[ \int_{\R^2} u(x) H(x)\ dx \leq \frac\alpha{2}\int_{\R^2}u^2(x)\ dx + \frac1{2\alpha}\int_{\R^2}H^2(x)\ dx\ .  \]

Let us mention to conclude this Appendix that there exists a third strategy to prove Lemma \ref{lem:Besselinequality}, which is based on a "weak logarithmic HLS inequality" (see \cite[Theorem 3]{Beckner}).

\end{document}